\documentclass{amsart}

\usepackage{epsfig,pinlabel}

\newtheorem{thm}{Theorem}[section]
\newtheorem{lemma}[thm]{Lemma}
\newtheorem{prop}[thm]{Proposition}

\theoremstyle{definition}\newtheorem*{ex}{Example}

\newcounter{Tlistc}
\newenvironment{Tlist}
    	{\setcounter{Tlistc}{0}
	 \begin{list}{(T\arabic{Tlistc})}
	{\usecounter{Tlistc}}}{\end{list}}
	 
\newenvironment{romanlist}
	{\begin{enumerate}
	}
	{\end{enumerate}}

\newcounter{ticklistc}
\newenvironment{ticklist}
    	{\setcounter{ticklistc}{0}
	 \begin{list}{-}
	{\usecounter{ticklistc}}}{\end{list}}

\newcommand{\Zz}{\mathbb Z}
\newcommand{\sur}{\Sigma}
\newcommand{\Int}{\text{Int}}
\newcommand{\Ker}{\text{Ker}}
\newcommand{\Rr}{\mathbb R}
\newcommand{\Cc}{\mathbb C}
\newcommand{\eps}{\varepsilon}
\newcommand{\sgn}{\mathrm{sgn}}
\newcommand{\rk}{\mathrm{rank}\,}
\newcommand{\Hom}{\mathrm{Hom}}
\newcommand{\spa}{\mathrm{span}\,}

\begin{document}

\title{A generalization of several classical invariants of links}

\author{David Cimasoni}
\address{David Cimasoni, Department of Mathematics, UC Berkeley, 970 Evans Hall, Berkeley, CA 94720, USA}
\email{cimasoni@math.berkeley.edu}
\thanks{The first author is supported by the Swiss National Science Foundation.}

\author{Vladimir Turaev}
\address{Vladimir Turaev, IRMA, CNRS et Universit\'e Louis Pasteur, 7 rue Ren\'e Descartes, 67084 Strasbourg Cedex, France}
\email{turaev@math.u-strasbg.fr}

\date{\today}


\begin{abstract}
We extend several classical invariants of links in the 3-sphere to links in so-called quasi-cylinders.
These invariants include the linking number, the Seifert form, the Alexander module, the Alexander-Conway polynomial
and the Murasugi-Tristram-Levine signatures.
\end{abstract}

\maketitle

\section*{Introduction}

The aim of this paper is to introduce a generalization of several classical knot 
and link invariants including the linking numbers, the Alexander-Conway polynomial, and 
the Murasugi-Tristram-Levine signatures of links in Euclidean 3-space. These invariants are 
generalized to links in so-called quasi-cylindres. A quasi-cylinder over a commutative ring $R$ is an 
oriented 3-manifold $M$ endowed with a submodule $V$ of the $R$-module $H_1(\partial M;R)$ such that the 
inclusion homomorphism $V\to H_1(M;R)$ is an isomorphism. The main example is the cylinder 
$M=\sur\times [0,1]$ where $\sur$ is an oriented surface and $V=H_1(\sur\times 0;\Zz)\subset H_1(\partial M;\Zz)$.
(Here, $R=\Zz$). For homologically trivial links in a quasi-cylinder $M$ 
with $H_2(M;\Zz)=0$,  we define a generalized 
Seifert form and derive from it several other invariants, namely,  a generalized Alexander-Conway 
polynomial  and  generalized  signatures. The most interesting feature of our 
invariants is the appearance of additional parameters which are absent in the classical case. 

The organization of the paper is as follows. In Section \ref{section:lk} we introduce 
generalized linking numbers of links in quasi-cylinders. In Section \ref{section:forms} we define the 
generalized Seifert form for links in quasi-cylinders. In Section \ref{section:Alex} we study the 
derived Alexander invariants. In Section \ref{section:genus} we discuss simple estimates of the 
link genus. In Section \ref{section:concordance} we study the  concordance of links.
In Section \ref{section:signature} we  consider the signatures. In Section \ref{section:colored} we introduce a 
multivariable extension of the theory. In Section \ref{section:gene} we discuss various generalizations of our invariants and in particular 
an extention to homologically non-trivial links. 
 
In this paper, all manifolds are smooth. The boundary of an oriented manifold 
is oriented {\em via\/} the \lq\lq outward normal vector first" convention. 

Throughout the paper, we fix a commutative ring with unity $R$.

\section{Linking numbers and quasi-cylinders}\label{section:lk}

\subsection{Knots and links}
By a {\em link\/} in an oriented 3-manifold $M$, we mean a finite system of 
disjoint  oriented circles embedded in $\Int(M)=M-\partial M$. Each link $L$ 
in $M$ viewed as  a geometric 1-cycle represents a homology class $[L]\in H_1(M;R)$.
A  link $L$ is {\em $R$-homologically trivial\/} if $[L]=0$. For $R=\Zz$, we 
say simply {\em homologically trivial}.

Two links $L$ and $L'$ in $M$ are said to be {\em ambient isotopic\/} if there is an ambient
isotopy $h_t\;(0\le t\le 1)$ of $M$, keeping $\partial M$ fixed, such that $h_0=id$, $h_1(L)=L'$,
and $h_1\vert_L\colon L\cong L'$ is orientation-preserving.

A {\em knot\/} is a link consisting of a single circle. Let us stress that 
all knots and links in this paper are oriented.

\subsection{Linking numbers}
The classical linking number of disjoint $R$-homologically 
trivial knots $K,L$ in an oriented 3-manifold $M$  is defined by
$$
lk (K, L)= K \cdot B =B \cdot K\in R
$$
where $\cdot=\cdot_M$ is the standard homological intersection in $M$ and $B $ 
is a  2-chain in $M$ (with coefficients in $R$) such that $\partial B = L$. The 
independence of the choice 
of $B$ follows from the fact that given another 2-chain $B'$ with $\partial 
B'=L$, 
one has $K \cdot B-K \cdot B'=K \cdot b$ where $b= B -B' $ is a 2-cycle in $M$. 
The 
$R$-homological triviality of $K$ implies that $K \cdot b=[K]\cdot b=0$.
One easily checks the symmetry $lk (K, L)=lk (L, K)$.

We introduce a {\em generalized linking number\/} as follows. Suppose   that 
$\partial M\neq \emptyset$ and   denote by $c$ the inclusion homomorphism 
$H_1(\partial 
M;R)\to H_1(  M ;R)$. Fix a submodule $V$ of the
$R$-module $ H_1(\partial M;R)$ such that   $V\cap \Ker (c) =0$.  For    disjoint     
knots $K, L$ in   $M$ such that $[K],[L] \in c(V)$, set
$$
lk_V(K,L)=K\cdot B=B\cdot K\in R
$$
where $B $ is any 2-chain in $M$ (with coefficients in $R$) such that $\partial B=L-v$
for a 1-cycle $v$ on $\partial M$ representing an element of $V$. The 
homological intersection   $K \cdot B$ does not depend on the choice of $B$. 
Indeed, 
consider  another 
2-chain $B'$ in $M$ with $\partial B' = L -v'$ where $v'$ is a 1-cycle on 
$\partial 
M$ representing an element of $V$. Then $b=B-B'$ is a relative 2-cycle in $(M, 
\partial M)$ in the sense that its boundary lies on $\partial M$.  Let $u$ be a 
1-cycle on $\partial M$ whose homology class $[u]\in H_1(\partial M;R)$ satisfies   
$c([u])=[K]$ and let  $\tilde u$ be a 1-cycle   obtained    by pushing $u$ 
slightly 
inside $\Int(M)$. Then
$$
K\cdot B-K \cdot B'=K \cdot b=\tilde u\cdot b= u\cdot_{\partial M}\partial b
$$ 
where $\cdot_{\partial M}$ is the homological intersection of 1-cycles in 
$\partial M$. We have $u \cdot_{\partial M} \partial b =0$ since $[\partial b]=[v-v']\in V\cap\Ker (c)=0$.  
 
The linking number $lk_V$   satisfies
$$
lk_V (L, K)=lk_V (K, L)+u \cdot_{\partial M} v.
$$
where $u, v$ are 1-cycles on $\partial M$ representing elements of $V$ 
homological to 
$K, L$ respectively. Indeed, let  $\tilde u$ be a 1-cycle in $\Int(M)$  obtained  
from $u$ as above and let    $A $ be a  2-chain in $M$ with $\partial A = K-\tilde u$. Then $A$ is disjoint form $v$ and therefore 
\begin{eqnarray*}
lk_V (L, K)&=&L \cdot A=(L -v ) \cdot A  = \partial B\cdot A  =B\cdot\partial A\\
&=&B\cdot K- B\cdot \tilde u =B \cdot K - v \cdot_{\partial M} u\\
&=&lk_V (K, L)+u \cdot_{\partial M} v.
\end{eqnarray*}
It is clear that $lk_V(K,L)$ is invariant under deformations of $K$ and $L$ in 
$M$ 
keeping them disjoint. If $K,L$ are $R$-homologically trivial (in particular, if 
they lie in a 
3-ball inside $M$), then $lk_V(K,L)=lk(K,L)$. 

The definition of $lk_V(K,L)$ extends in the obvious way to the case where $K,L$ 
are disjoint 1-cycles in $M$.

\subsection{Quasi-cylinders}
By a {\em quasi-cylinder\/} (over $R$), we mean a pair 
consisting of a connected  oriented 3-manifold $M$ with non-empty boundary and a 
submodule $V$ of the 
$R$-module $H_1(\partial M; R)$ such that the restriction of the inclusion homomorphism
$H_1(\partial M;R)\to H_1(M;R)$ to $V$ yields an isomorphism $V\to H_1( M;R)$.
The inverse isomorphism is denoted $d_V$. 

The constructions of the previous section show that for a quasi-cylinder  
$(M,V)$ over $R$ and any 
disjoint knots $K,L$ in $M$, we have a well-defined linking number $lk_V(K,L)\in 
R$ satisfying
$$
lk_V (L, K)=lk_V (K, L)+d_V([K]) \cdot_{\partial M} d_V([L]).
$$
 
We say that a quasi-cylinder $(M,V)$ has trivial 2-homology if $H_2(M)=0$. 
Here and below, the unspecified group of coefficients in homology/cohomology is $\Zz$.

Note the following lemma.

\begin{lemma}\label{lemma:homology}
Let $(M, V)$ be a quasi-cylinder such that $M$ is 
compact. The equality $H_2(M)=0$ holds if and only if $\partial M$ is 
connected.
\end{lemma}
\begin{proof}
The components of $\partial M$ represent 
elements of $H_2(M)$  subject to only one relation: their sum   is equal to zero. 
Therefore the 
equality $H_2(M)=0$ implies that  
  $\partial M$ is connected. Let us prove the converse. Since the inclusion 
homomorphism $ 
H_1(\partial M;R)\to H_1(  M;R)$ is onto and the inclusion homomorphism $ 
H_0(\partial 
M;R)\to H_0(  M;R)$ is an isomorphism, the homology sequence of the pair $(M, 
\partial M)$ 
gives that $H_1(M,\partial M;R)=0$. Observe that $H_1(M,\partial 
M;R)=R\otimes_{\Zz} 
H_1(M,\partial M)$. Therefore the group $H_1(M,\partial M)$ is finite and 
$$
H_2(M)=H^1 (M, \partial M)=\Hom (H_1(M,\partial M), \Zz)=0.
$$
\end{proof}

\subsection{Examples}
1. The pair consising of a 3-ball $D^3$ and   
$V= H_1(\partial D^3)=0$ is a quasi-cylinder over $ \Zz$. Clearly,  
$lk_V (K, L)=lk  (K, L)\in \Zz$ is the usual linking number of knots in the 
3-ball.  

2. Let $\sur$ be a connected oriented surface and $M=\sur \times [0,1]$ with 
product orientation. Set $V=H_1( \sur \times 0;R) \subset H_1(\partial M;R)$. It is clear that $(M,V)$ is 
a quasi-cylinder. The linking number $lk_V (K, L)$ of knots $K,L\subset M$ can be 
computed as follows. Present the link $K\cup L$ by a link diagram on $\sur$. Let $k, l$ be 
the components of the diagram representing $K$ and $L$, respectively. Then $lk_V (K, L)=n_+-n_-$ 
where $n_+$ (resp.\ $n_-$) is the number of positive (resp.\ negative) crossing points on the 
diagram where $k$ passes under  $l$. The quasi-cylinder $(\sur \times [0,1], V)$ has 
trivial 2-homology if and only if $\partial \sur\neq \emptyset$.

3. Let $N$ be an $R$-homology 3-sphere, i.e.,     a closed oriented 3-manifold 
with $H_\ast 
(N; R)=H_\ast (S^3; R)$. Let $G$ be a non-empty finite   graph in   $N$ and  
$M\subset N$  be its exterior, that is the complement of an open regular neighborhood of 
$G$ in $N$. Let $V\subset H_1(\partial M; R)$ be the submodule generated by the homology 
classes of the meridians of the edges of $G$. Then the pair $(M,V)$ is a quasi-cylinder. For any
knots $K,L\subset M$, we have $lk_V (K, L)=lk  (K, L)\in R$ where the right-hand side is the 
linking number of $K,L$ in $N$. The quasi-cylinder $(M,V)$ has trivial 2-homology if and only if $G$ is connected.

\subsection{Remark}
The constructions above suggest that the definition of the Milnor numbers 
of classical links may be  extended to links in quasi-cylinders. We shall not 
pursue this line here.

\section{Generalized Seifert forms}\label{section:forms}

\subsection{Bilinear forms associated with surfaces}
Let $(M,V)$ be a quasi-cylinder over $R$. Consider an oriented embedded surface $F\subset\Int(M)$ (possibly,
$\partial F\neq\emptyset$). For a 1-cycle $a$ on $F$, denote by $a^+$ (resp.\ $a^-$) the 
1-cycle in $\Int(M)\backslash F$ obtained by pushing $a$ along the positive  
(resp.\ negative) normal direction on $F$ in $M$.  
For 1-cycles $a,b$ on $F$ representing homology classes $[a], [b] \in H_1(F;R)$, set 
$$
\vartheta ([a],[b])=lk_V(a^+,b).
$$
This number only depends on the homology classes of $a$ and $b$ in $H_1(F;R)$. Indeed, if $a_1,a_2,b_1$ and $b_2$
are $1$-cycles on $F$ such that $a_1-a_2=\partial A$ and $b_1-b_2=\partial B$ for some $2$-chains $A$ and $B$ in $F$, then
\begin{eqnarray*}
lk_V(a_1^+,b_1)-lk_V(a_2^+,b_2)&=&lk_V(a_1^+-a_2^+,b_1)+lk_V(a_2^+,b_1-b_2)\\
&=&\partial A^+\cdot B_1+a_2^+\cdot B=A^+\cdot\partial B_1+a_2^+\cdot B\\
&=&A^+\cdot b_1-A^+\cdot v_1 +a_2^+\cdot B,
\end{eqnarray*}
where $A^+$ denotes the $2$-chain $A$ pushed along the positive normal direction on $F$ in $M$.
Since $A^+,a_2^+\subset \Int(M)\setminus F$ and $b_1,B\subset F$, $v_1\subset\partial M$, these three intersection number are zero.
Hence, we have a  well-defined bilinear form 
$$
\vartheta =\vartheta_F\colon H_1(F;R) \times H_1(F;R)\to R.
$$
We call it the {\em generalized Seifert form\/} of $F$.  

\begin{lemma}\label{lemma:relation}
Let $ d \colon H_1(F;R)\to V$ be the composition 
of the 
inclusion homomorphism $H_1(F;R)\to H_1(M;R)$ with the isomorphism 
$d_V\colon H_1(M;R)\to V$.
For all $a,b\in H_1(F;R)$,
$$
a\cdot_F b=\vartheta (a,b)-\vartheta (b,a)-d(a) \cdot_{\partial M} d(b).
$$
\end{lemma}
\begin{proof} By abuse of notation, we shall denote   1-cycles    
representing $a,b, d(a), d(b)$ by the same symbols $a,b, d(a), d(b)$. Consider 
the 
bilinear form $\vartheta^-\colon H_1(F;R) \times H_1(F;R)\to R$  defined as 
$\vartheta^+=\vartheta$ 
but 
using  $a^-$ instead of  $a^+$. We claim that
$$
\vartheta^+(a,b)-\vartheta^-(a,b)=a \cdot_F b. \eqno(2.a)
$$
Indeed, let $B$ be a $2$-cycle in $M$ such 
that $\partial B=b-d(b)$, and let $\alpha$ be the $2$-cycle
$[-1,1]\times a$ in $\Int(M)$ with  
$\partial\alpha=a^+-a^-$. We have
\begin{eqnarray*}
\vartheta^+(a,b)-\vartheta^-(a,b)&=& a^+\cdot B - a^- \cdot B = \partial\alpha 
\cdot B \cr
&=& \alpha\cdot \partial 
B = \alpha\cdot b - \alpha \cdot d(b).
\end{eqnarray*}
Now, $\alpha$ and $d(b)$ are disjoint, so $ \alpha \cdot d(b)=0$. 
Since $ \alpha\cdot b= a\cdot_F b$, this gives (2.a).

We now verify that 
$$
\vartheta^+(a,b)-\vartheta^-(b,a)=d(a)\cdot_{\partial M} d(b). \eqno(2.b)
$$
Indeed,
\begin{eqnarray*}
\vartheta^+(a,b)&=&lk_V(a^+,b)=lk_V(b, a^+)+ d(a) \cdot_{\partial M} d(b)\\
&=&lk_V(b^-, a)+ d(a) \cdot_{\partial M} d(b)=\vartheta^-(b,a)+d(a)\cdot_{\partial M} d(b).
\end{eqnarray*}
Combining formulas (2.a) and  (2.b), we obtain the claim of the lemma.
\end{proof} 
 
\subsection{Algebraic digression}
Let $W$ be an arbitrary $R$-module. A {\em Seifert triple\/} over $W$ is a triple 
$(H,\vartheta,d)$, where $H$ is a free $R$-module of finite rank, $\vartheta$ a  bilinear form  
$H\times H \to R$, and $d$ a homomorphism $ H\to W$. Two Seifert triples $(H_1,\vartheta_1,d_1)$,
$(H_2,\vartheta_2,d_2)$ over $W$  are isomorphic if there is an 
$R$-isomorphism $f\colon H_1 \to H_2$ such that $\vartheta_2 \circ (f \times f)= \vartheta_1$ and $d_2\circ f=d_1$.

A Seifert triple $(H',\vartheta',d')$ is obtained from a  Seifert 
triple $(H,\vartheta,d)$ by an {\em elementary enlargement\/} (and
$(H,\vartheta,d)$ from $(H',\vartheta',d')$ by an {\em elementary reduction\/}) 
if the following conditions hold:
$H'= H\oplus R  a\oplus Rb$, $d'|_H=d$, $d'(b)=0$, $\vartheta'|_{H\times H}=\vartheta$,
$\vartheta'(H,b)= \vartheta'(b, H)=\vartheta'(b, b)=0$ and either 
$\vartheta'(a,b)=1,\vartheta'(b, a)=0$ or $\vartheta'(a,b)=0,  \vartheta'(b, a)=1$. If $h$ is a   
basis of $H$, then $h \cup \{a,b\}$ is a basis of $H'$ and the 
matrix $\Theta'$ of $\vartheta'$ with respect to $h \cup \{a,b\}$ is computed from the matrix 
$\Theta$ of $\vartheta$ with respect to   $h$  by
$$
\Theta'=\begin{pmatrix}\Theta&\star&0\cr\star&\star&0\cr 
0&1&0\end{pmatrix}
\quad\hbox{or}\quad
\Theta'=\begin{pmatrix}\Theta&\star&0\cr\star&\star&1\cr 
0&0&0\end{pmatrix}.
$$
We say that two Seifert triples over $W$ are 
{\em equivalent\/} if they can be related by a finite sequence of isomorphisms 
and elementary enlargements and reductions.

\subsection{Surgeries on surfaces}
Given a quasi-cylinder  $(M, V)$ over $R$ and a compact connected oriented    
surface  $F\subset \Int(M)$, the constructions above yield a Seifert triple $(H_1 (F;R),\vartheta,d)$ over $V$.
Note that the $R$-module $H_1(F;R)$ is free of rank $2g+m-1$, where $g$ is the
genus of $F$ and $m$ is the number of connected components of $\partial F$.
Suppose that a surface $F'\subset\Int(M)$ is obtained from $F$ by surgery along an embedded arc in $\Int(M)$ meeting $F$ exactly at 
its endpoints and approaching $F$ either from the positive side or from the 
negative side at both endpoints. The transformation $F\mapsto F'$ and the inverse 
transformation are called {\em surgeries}. It is easy to see that  the Seifert 
triple of  $F'$ is obtained from the Seifert triple of $F$ by an elementary enlargement. 
Therefore the equivalence class of the Seifert triple of an embedded 
surface is invariant under surgeries.

Observe that any given class in $H_2(M)$ can be realized by a closed connected 
oriented surface and any two such surfaces are related  by a sequence  of   
surgeries and isotopies. This leads to algebraic invariants of integral 2-homology 
classes of quasi-cylinders. We shall however focus on quasi-cylinders with trivial 
2-homology. 

\subsection{Seifert forms of links}
A {\em Seifert surface\/} for a link $L$ in a 3-manifold $ M $ is a compact 
connected oriented surface in $\Int(M)$ that has $L$ as its oriented
boundary. Clearly, if $L$ has a Seifert surface, then $[L]=0$ in $H_1(M)$. It is well-known that   
this is the only obstruction. For completeness, we outline a proof.

\begin{lemma}\label{prop:surface}
Any homologically trivial link  in an oriented 3-manifold has a Seifert surface.
\end{lemma}
\begin{proof} Let $L$ be a homologically trivial link in  an oriented 3-manifold $M$. Then 
$L$ is homologically trivial in a compact 3-dimensional submanifold $M'$ of $M$ such that 
$M'\supset L$. Let $N $ be a closed tubular neighborhood of $L$ in $\Int(M')$. 
Set $X=M'\setminus\Int(N)$. Since $[L]=0\in H_1(M')$,
an appropriate choice of longitudes of $L$ gives a link $L'\subset\partial N\subset\partial X $ 
whose  class in $H_1(X)$ is equal to $0$. 
Then $[L'] \in H_1( \partial X)$ is the boundary of an element of 
$ H_2(X,\partial X)=H^1(X)$. The latter   is the pull-back of a generator of 
$H^1(S^1 )=\Zz$ under a  map $X\to S^1$. For an appropriate choice of this map, 
the preimage of a point of $S^1$   is  a  compact oriented surface bounded by $L'$ in $X$.
Adding if necessary 1-handles to this surface one can 
make it connected. The resulting surface extends to a Seifert surface for $L$ in $M'$.
\end{proof}

Given two Seifert surfaces $F, F'$ for a link $L$ in an oriented 3-manifold $M$, 
the  union $F\cup (-F')$ is a closed oriented
surface representing an element of $H_2(M)$. This element is an obstruction to 
transforming $F$ into $F'$ by surgeries. It is well-known that this is the only 
obstruction (see e.g.\ \cite[p.\ 64]{Kaw}). In particular, if $H_2(M)=0$, then 
$F,F'$   
can be related by a finite sequence of surgeries  and ambient 
isotopies  in $M$
(which can be chosen to keep $\partial M$ fixed). Combining this fact  with the 
observations above, we obtain the following.

\begin{thm}\label{thm1} Let $(M, V)$ be a quasi-cylinder over   
$R$ with $H_2(M)=0$. For any homologically trivial link $L\subset M$, the   
equivalence class of the Seifert triple of a Seifert surface for $L$  does not depend on the 
choice of the surface and provides an isotopy invariant of $L$. 
\end{thm}

\section{Alexander invariants}\label{section:Alex}

Throughout this section, $(M,V)$ is a quasi-cylinder over $R$ with $H_2(M)=0$. 

\subsection{The Alexander module}\label{sec2.1}
Fix a commutative unital ring $R'$ containing $R$ as a subring. We also fix an $R$-bilinear pairing $\psi\colon V \times V \to R'$.  Consider a 
homologically trivial link $L$ in $M$. Let  $(H,\vartheta\colon H\times H\to R,d\colon H \to V)$ be 
the Seifert triple associated with a Seifert surface for $L$. Let $\Theta$ and $\Psi$ be the  
matrices of the bilinear forms $\vartheta$ and $\psi\circ(d\times d)$ with 
respect to a basis of $H$. The {\em Alexander module\/} $\mathcal{A}_\psi(L)$ of $L$ is the 
$R'[t,t^{-1} ]$-module  presented by the matrix $t\Theta-\Theta^T+ \Psi$, where the superscript $T$ denotes the  matrix transposition.

\begin{prop}\label{prop:module}
The Alexander module is an isotopy invariant of $L$.
\end{prop}
\begin{proof}
Obviously, this module does not depend on the choice of a basis of $H$. 
By Theorem \ref{thm1}, we just need to check that if 
$(H',\vartheta',d')$
is obtained from $(H,\vartheta,d)$ by an elementary enlargement, then 
the corresponding matrices $\Gamma'=t\Theta'-(\Theta')^T+ \Psi'$ and
$\Gamma=t\Theta-\Theta^T+ \Psi$
present isomorphic $R'[t,t^{-1} ]$-modules. Clearly,
$$
\Psi'=\begin{pmatrix}\Psi&\star&0\cr\star&\star&0\cr 
0&0&0\end{pmatrix}.
$$
Therefore,
$$
\Gamma'=\begin{pmatrix}\Gamma&\star&0\cr\star&\star&-1\cr 0&t&0\end{pmatrix}
\quad\hbox{or}\quad
\Gamma'=\begin{pmatrix}\Gamma&\star&0\cr\star&\star&t\cr 0&-1&0\end{pmatrix}.
$$
In both cases, the corresponding modules over $R'[t,t^{-1} ]$ are 
isomorphic.
\end{proof}

For $M=D^3$, $V=0$, $R'=R=\Zz$ and  $\psi=0$,  the module $\mathcal{A}_\psi(L)$ is the usual Alexander module.
  
Mimicking the standard definitions, we can introduce the Alexander ideals and Alexander polynomials of $L$
(provided $R'$ is a unique factorization domain). In particular, the first Alexander polynomial of $L$ can be defined as the 
determinant of a square presentation matrix of  $\mathcal{A}_\psi(L)$. This polynomial is an element of the ring $R'[t,t^{-1}]$
defined up to multiplication by units of this ring. As in the classical case, the first Alexander polynomial has a canonical
normalization which we now discuss. 

\subsection{The Alexander-Conway polynomial}\label{sec2.2}
Using the notation of the previous subsection, we define the {\em (extended) Alexander-Conway polynomial\/} of $L$ by
$$
\Delta_{ L,\psi}(t)=\det(t^{1/2}\Theta-t^{-1/2}\Theta^T+ t^{-1/2}\Psi).
$$
As in the proof of Proposition \ref{prop:module}, one checks that this element of $R'[t^{1/2},t^{-1/2}]$ is a well-defined isotopy invariant of $L$.

Observe that the size of the matrices $\Theta, \Psi$ is equal to $2g +m-1$, where 
$g$ is the genus of the Seifert surface and $m$ is the number of components of $L$. Therefore  
$$
\Delta_{ L,\psi}(t )=t^{\frac{1-m}{2} -g }\det(t\Theta-\Theta^T+ \Psi).
$$
Thus, $\Delta_{L,\psi}(t)\in R'[t,t^{-1}]$ for odd $m$ and $t^{1/2}\Delta_{L,\psi}(t)\in R'[t,t^{-1}]$ for even $m$.

We now establish a skein formula for $\Delta_{L,\psi}(t)$.

\begin{prop}\label{Conway}
Let $L_+$, $L_-$ and $L_0$ be homologically trivial links in $M$ which coincide 
everywhere except in a small $3$-ball where they are related as illustrated below.
\begin{figure}[h]
\labellist\small\hair 2.5pt
\pinlabel {$L_+$} at 87 -20
\pinlabel {$L_-$} at 340 -20
\pinlabel {$L_0$} at 570 -20
\endlabellist
\centerline{\psfig{file=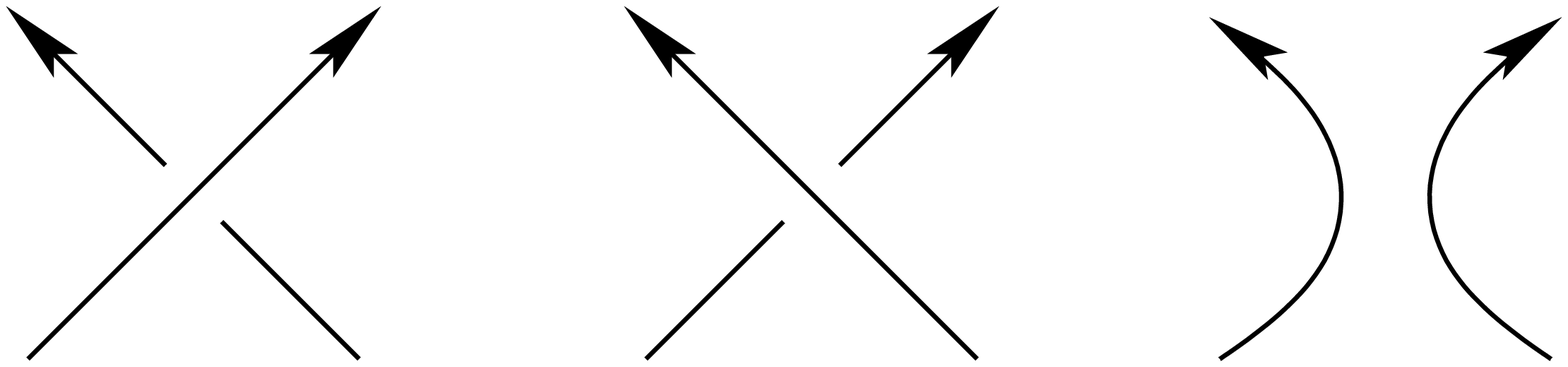,height=1.5cm}}
\end{figure}

\noindent Then, the corresponding Alexander-Conway polynomials satisfy the following relation:
$$
\Delta_{L_+,\psi}(t)-\Delta_{L_-,\psi}(t)=(t^{-1/2}-t^{1/2}) \, \Delta_{L_0,\psi}(t,s).
$$
\end{prop}
\begin{proof}
Let $F_0$ be a Seifert surface for $L_0$. Then a Seifert surface $F_+$ for $L_+$ (resp.\ $F_-$ for $L_-$)
is obtained from $F_0$ by adding a band in the small $3$-ball with one negative (resp.\ positive) half-twist. Since $F_0$
is connected, a basis for $H_1(F_+;R)$ (resp.\ for $H_1(F_-;R)$) is obtained 
from a basis for $H_1(F_0;R)$ by adding a $1$-cycle $a_+$ (resp.\ $a_-$).
Clearly, $a_+$ and $a_-$ can be chosen to coincide as $1$-cycles in 
$M$. Let $v$ be a 1-cycle on $\partial M$ with coefficients in $R$ such that $[v]\in V\subset H_1(\partial M;R)$
and $v$ is homologous to  $a_+=a_-$ in $M$. Let $B$ be a $2$-cycle in $M$ such that $\partial B=a_{\pm}- v$. Then
$$
\vartheta_{F_+}(a_+,a_+)-\vartheta_{F_-}(a_-,a_-)= a_+^+\cdot B- a_-^+\cdot B=(a_+^+-a_-^+)\cdot B =-1.
$$
This leads to the following equalities between the corresponding matrices:
$$
\Theta_{F_+}=\begin{pmatrix}\Theta_{F_0}&v\cr 
w &\alpha\end{pmatrix},\;\;\Theta_{F_-}=\begin{pmatrix}\Theta_{F_0}&v\cr 
w &\alpha+1\end{pmatrix},
\;\;\hbox{and}\quad 
\Psi_{F_+}=\Psi_{F_-}=\begin{pmatrix}\Psi_{F_0}&x\cr y&\beta\end{pmatrix}
$$
for some   $\alpha\in R,\beta\in R'$, column  $v $ and row $w$ over $R$, and column $x$ and row $y$ over $R'$. The skein 
formula follows.
\end{proof}

The skein formula implies in particular that $\Delta_{L,\psi}(1)\in R'$ is unchanged when
one replaces an undercrossing by an overcrossing. Hence, it depends only on the homotopy type 
of the components of $L$.

If $L'$ is a link in an oriented 3-ball $D^3$ and $L $ is the image of $L'$ 
under an orientation preserving embedding $D^3\hookrightarrow M$, then 
$\Delta_{L,\psi}(t )=\Delta_{L'}(t)$ is the usual Conway-normalized Alexander polynomial of $L'$.  

\subsection{A special case}
Let $R'=R[s_1,\dots,s_n]$ be the polynomial ring over $R$ generated by $n$  commuting variables $s_1,\dots, s_n$. Let
$\psi_1,\dots, \psi_n\colon V\times V \to R$ be  bilinear forms. We can apply the definitions and results of the previous subsections to the
bilinear form
$$
\psi=s_1\psi_1+ \cdots + s_n \psi_n \colon  V\times V \to R'.
$$
This gives a polynomial invariant
$$
\Delta_{L,\psi_1,\dots,\psi_n}(t,s_1,\dots,s_n)=\Delta_{L,\psi}(t)=\det(t^{-1/2}\Theta-t^{1/2}\Theta^T+t^{-1/2}\sum_{i=1}^n s_i\Psi_i),
$$
where $\Theta$ and $\Psi_i$ are the matrices of the bilinear forms $\vartheta$ and $\psi_i\circ(d\times d)$ with 
respect to a basis of $H$. The polynomial $\Delta_{L,\psi_1,\dots,\psi_n}(t,s_1,\dots,s_n)$
lies in $R[t,t^{-1},s_1,\dots,s_n]$ for odd $m$ and in $t^{1/2}\times R[t,t^{-1},s_1,\dots,s_n]$ for even $m$. 

The degree in $s_i$ of $\Delta_{ L,\psi_1,\dots, \psi_n}(t,s_1,\dots, s_n )$ is bounded from above by a number independent of $L$.
Namely, this degree is smaller than or equal to
the rank of the form $\psi_i$. Indeed, it follows from the definitions that
$$
\deg_{s_i}\Delta_{L,\psi_1,\dots, \psi_n}(t,s_1,\dots,s_n)\leq \rk(\Psi_i)=\rk \psi_i.
$$

For a bilinear form $\psi\colon V \times V\to R'$ we denote by $\psi^T$ its transpose defined by $\psi^T(a,b)=\psi(b,a)$ for $a,b\in V$.

\begin{prop}\label{prop:properties}
$$
\Delta_{L,\psi_1,\dots,\psi_n}(t^{-1},s_1t^{-1},\dots,s_n t^{-1})=(-1)^{m-1}\Delta_{L,-\psi^T_1,\dots,-\psi^T_n}(t,s_1,\dots,s_n).
$$
\end{prop}
\begin{proof}
Transposing matrices, we obtain 
\begin{eqnarray*}
\Delta_{ L,\psi_1,\dots, \psi_n}(t^{-1},s_1t^{-1},\dots, s_n t^{-1})&=&\det(t^{-1/2}\Theta-t^{1/2}\Theta^T+t^{1/2}\sum_{i=1}^n s_it^{-1}\Psi_i)\\
&=&\det(t^{-1/2}\Theta^T-t^{1/2}\Theta+t^{-1/2} \sum_{i=1}^n s_i \Psi_i^T) \\
&=&(-1)^{m-1}\Delta_{L,-\psi^T_1,\dots,-\psi^T_n}(t,s_1,\dots,s_n).
\end{eqnarray*}
\end{proof}
For example, if $\psi_i$ is symmetric for $i=1,\dots, p$ and skew-symmetric for $i=p+1,\dots, n$, then    
$$
\Delta_{L,\psi_1,\dots,\psi_n}(t^{-1},s_1t^{-1},\dots,s_n t^{-1})=(-1)^{m-1}\Delta_{L,-\psi_1,\dots,-\psi_p,\psi_{p+1},\dots,\psi_n}(t,s_1,\dots,s_n)
$$
$$
=(-1)^{m-1}\Delta_{L,\psi_1,\dots,\psi_n}(t,-s_1,\dots,-s_p,s_{p+1},\dots,s_n).
$$

The polynomial  $\Delta_{L,\psi_1,\dots,\psi_n}(t,s_1,\dots,s_n)$   
leads to other polynomial invariants of $L$. First of all, we can expand
$$
\Delta_{L,\psi_1,\dots,\psi_n}(t,s_1,\dots,s_n)=\sum_{i_1,\dots,i_n\geq 0}\Delta_L^{(i_1,\dots, i_n)}(t)\,s_1^{i_1}\cdots s_n^{i_n},
$$
where $\Delta_L^{(i_1,\dots, i_n)} (t)$ belongs to $R[t,t^{-1}]$ for odd $m$ and to $t^{1/2}\times R[t,t^{-1}]$ for even $m$.
The sum on the right-hand side is finite since 
$\Delta_L^{ (i_1,\dots, i_n)} (t)=0$ provided $i_k> \rk \psi_{i_k}$ for some $k=1,\dots,n$. For any triple $(L_+, L_-, L_0)$ as in Proposition 
\ref{Conway} and for any $i_1,\dots, i_n\geq 0$,
$$
\Delta_{L_+}^{(i_1,\dots,i_n)}(t)-\Delta_{L_-}^{(i_1,\dots,i_n)}(t)=(t^{1/2}-t^{-1/2})\,\Delta_{L_0}^{(i_1,\dots,i_n)}(t).
$$

Another interesting restriction of $\Delta_{L,\psi_1,\dots,\psi_n}$ is obtained by the substitution $t=1$.
By the skein relation, the resulting polynomial depends only on the homotopy type of the components of $L$.

If $V$ is a free module, then we can take as $\psi_1,\dots,\psi_n$ 
a basis in the $R$-module of bilinear forms $V\otimes V\to R$. This results
in a link polynomial on $1+v^2$ variables, where $v$ is the rank of $V$. 

\begin{ex}
Take $n=1$ and let $\psi_1=\cdot_{\partial M}\colon V \times V\to  R$ be the homological intersection on $\partial M$ restricted to $V$. This gives a polynomial invariant  $\Delta_{ L }(t,s)= \det(t^{1/2}\Theta-t^{-1/2}\Theta^T+ t^{-1/2}s \Psi )$ 
where $\Theta$ and $\Psi$ are the  
matrices of the bilinear forms $\vartheta$ and $\cdot_{\partial M}\circ(d\times d)$ with 
respect to a basis of $H$. We leave to the reader to check the following three properties of $\Delta_{L}(t,s)$, where $m$
denotes the number of components of $L$:
\begin{ticklist}
\item{$\Delta_L(1,-1)=1$ if $m=1$ and $\Delta_L(1,-1)=0$ otherwise;}
\item{$\Delta_{ L }(t^{-1},s t^{-1} )=(-1)^{m-1}\Delta_{L}(t,s)$;}
\item{$\Delta_{-L}(t,s)=\Delta_L(t,-(s+t+1))$, where $-L$ is $L$ with opposite orientation.}
\end{ticklist}
We can sometimes explicitly compute $\Delta_L(t,s)$ for links $L$ represented by
simple closed curves on $\partial M$.  Let $\sur \subset \partial M$ be 
a compact  connected  surface  of genus g with boundary, such that the image of 
the inclusion homomorphism $H_1(\sur;R)\to H_1(\partial M;R)$ is 
contained in $V$. We endow $\sur$ with the orientation induced by the orientation 
on 
$\partial M$ (which in its turn is induced by the one on $M$). Let $\tilde  
{\sur}\subset 
\Int(M)$ be the oriented surface obtained by pushing $\sur$ inside $M$ and 
reversing its 
orientation.  Clearly,  $L=\partial  \tilde  {\sur}\subset M$ is a homologically 
trivial  
link with Seifert surface $\tilde  {\sur}$. 
It is easy to see that the form $\vartheta$ associated with  $\tilde  {\sur}$ 
is identically zero. It follows from the definitions that $\Delta_L(t,s)=t^{-g}s^{2g}$ if $L$ is a knot, and $\Delta_L(t,s)=0$ else.
\end{ex}

\subsection{Remark}
Let $\sur$ be a compact connected oriented surface of genus $g$ with 
$\partial \sur \neq \emptyset$. Consider the quasi-cylinder $M=\sur\times[0,1],V=H_1(\sur\times 0)$ over $\Zz$.
For any knot $K$ in $M$, the Laurent polynomial $\Delta=\Delta_K (t,s)\in\Zz[t,t^{-1},s]$ introduced in the previous example 
satisfies $\Delta(t^{-1},st^{-1})=\Delta(t,s)$, $\Delta(1,-1)=1$ and $\deg_s\Delta\leq 2g$. If $g=0$ (that is, if $\sur$ is a   
disc with holes), then these conditions characterize completely the polynomials $\Delta$ 
which can be realized as the Alexander-Conway polynomial of a
knot in $M$. Indeed, in this case $\Delta\in \Zz[t,t^{-1}]$, $\Delta(t^{-1})=\Delta(t),\Delta(1)=1$ so that $\Delta$ can be realized as the 
Alexander-Conway polynomial of a knot in a 3-ball in $M$. We do not know whether the conditions above are sufficient for $g>0$.

\section{The genus}\label{section:genus}

The genus of a homologically trivial link $L$ in a 3-manifold $M$ is defined by
$$
g(L)=\min\{\hbox{genus}(F)\,:\,\hbox{$F$ is a Seifert surface for $L$ in $M$}\}.
$$
If $M$ is a 3-ball, then Seifert proved that $g(L)\ge\frac{1}{2}(\spa\Delta_L(t)+1-m)$, where
$\spa$ is the usual span of a Laurent polynomial in one variable $t$ and $m$ is the number of components of $L$. This result
extends to our setting as follows.

\begin{prop}\label{prop:genus1}
Let $L$ be a homologically trivial $m$-component link in a 
quasi-cylinder $(M,V)$ with $H_2(M)=0$. Let $\psi\colon  V\times V \to R'$ be a pairing as in Section \ref{sec2.1}. Then 
$$
g(L)\ge\frac{1}{2}(\spa\Delta_{L,\psi}(t )+1-m).
$$
\end{prop}
\begin{proof}
Let $F$ be a Seifert surface for $L$ realizing the   genus 
$g(L) $, and let $\Theta,\Psi$ be corresponding matrices.
By definition of $\Delta_{L,\psi}(t )$,
\begin{eqnarray*}
\spa\Delta_{L,\psi}(t )&=&\spa\begin{vmatrix}t^{1/2}\Theta-t^{-1/2}\Theta^T+t^{-1/2}
\Psi\end
{vmatrix}\\
&=&\spa\begin{vmatrix}t\Theta-\Theta^T+ \Psi\end{vmatrix}\le\rk 
H_1(F)=2g(L)+m-1.
\end{eqnarray*}
The inequality follows.
\end{proof}

Consider now  a homologically trivial $m$-component link $L$ in $\sur\times 
[0,1]$, where $\sur $ is a compact connected oriented surface of genus $g$. The following  
algorithm (due to Seifert in the case where $\sur$ is a 2-disc) produces a Seifert 
surface for $L$ from a  connected diagram of $L$ on $\sur$. (A link  diagram is  connected if it cannot 
be presented as a   union of disjoint non-empty link diagrams.)  Let $n$ be the number of  
crossings on the   diagram. Smoothing these   crossings   in the unique way compatible with the 
orientation of $L$, one obtains a closed oriented 1-manifold  $\Gamma\subset\sur$ consisting of $\gamma\geq 1$ disjoint
simple closed curves on $\sur$. Note that $[\Gamma]=[L]=0\in H_1(\sur)$. 
Therefore, there is a finite collection of oriented connected subsurfaces $\Sigma_1,\dots,\Sigma_c$ of $\sur=\sur\times 0$ whose boundaries 
are disjoint and $\cup_i\partial\Sigma_i=\Gamma$. A Seifert surface $F$ for $L$ can be 
obtained from the $\Sigma_i$ by pushing their interiors into $\sur\times [0,1]$ and adding
a half-twisted band at each crossing.

\begin{prop}\label{prop:genus2}
Let $\gamma_0$ be the number of discs among the surfaces 
$\Sigma_1,\dots,\Sigma_c$. Then $\gamma_0\leq \gamma$ and 
$$
g(L)\le 1+\frac{1}{2}(n-\gamma-m)+(\gamma-\gamma_0)\max\{1,g \}.
$$
\end{prop}
\begin{proof} 
We have
$$
2-2g(F)-m=\chi(F)=\sum_{i=1}^c\chi(\Sigma_i)-n=2c-2\sum_{i= 1}^cg_i-\gamma-n
$$
where $g_i $ is the genus of $\Sigma_i$. Clearly $g_i\leq g$ and $g_i=0$ if $\Sigma_i $ is a disc. Hence
\begin{eqnarray*}
g(L)&\le&g(F)=1+\frac{1}{2}(n-\gamma-m)+(\gamma-c)+\sum_{i=1}^cg_i\\
&\le&1+\frac{1}{2}(n-\gamma-m)+(\gamma-c)+(c-\gamma_0) g.
\end{eqnarray*}
The inequalities $\gamma_0\le c\le \gamma$ now give the result.
\end{proof}

Note that if $g=0$, then $\gamma=\gamma_0$ and we obtain Seifert's inequality 
$g(L)\le 1+\frac{1}{2}(n-\gamma-m)$ for links in the 3-ball.

Combining Propositions \ref{prop:genus1} and \ref{prop:genus2}, we obtain in the case $\partial\Sigma\neq\emptyset$ that
$$
\spa\Delta_{L,\psi}(t)\le n+1-\gamma +2 (\gamma-\gamma_0)\max\{1,g \}.
$$

\section{Concordance invariants}\label{section:concordance}

Two links $L_0,L_1$ in a 3-manifold $M$ are {\em concordant\/} if there is a smooth oriented surface $S\subset M\times[0,1]$ such that 
$\partial S=(L_1\times 1)\cup(-L_0\times 0)$ and each component  of 
$S$ is an annulus with one boundary component on $M\times 0$ and the other one 
on $M\times 1$.  Cobordant links have the same number of components.

\begin{lemma}\label{lemma:conc}
Assume that $R$ is a principal ideal domain. Let $(M,V)$ be a quasi-cylinder over $R$ 
such that $M$ is compact and $H_2(M)=0$. Let $\psi\colon  V\times V \to R'$ be a bilinear pairing with values in an
integral domain $R'$ containing $R$ as a subring. Let $L_0,L_1$ 
be concordant homologically trivial links in $M$ and $F_0,F_1$ be their Seifert
surfaces with associated Seifert triples $(H_1(F_0;R ),\vartheta_0,d_0)$, 
$(H_1(F_1 ;R ),\vartheta_1,d_1)$. Then  there is a basis $x_1,\dots,x_{2g}$ of the $R$-module $H=H_1(F_0;R)\oplus H_1(F_1;R)$ such that the bilinear forms 
$$
\vartheta=(-\vartheta_0)\oplus\vartheta_1\quad {\text {and}} \quad \widetilde\psi= -(\psi\circ(d_0\times d_0))\oplus(\psi\circ(d_1\times d_1))
$$ 
satisfy $\vartheta(x_i,x_j)=\widetilde\psi(x_i,x_j)=0$ for all $i,j>g$.
\end{lemma}
\begin{proof}
Let $S\subset M\times [0,1]$ be a surface as in the definition of the link  
concordance.    Then 
$S\cup(F_0\times 0)\cup(-F_1\times 1)$ is a closed connected oriented surface in 
$M\times [0,1]$.

{\em Claim 1:\/} There is a compact   oriented $3$-manifold 
$N\subset M\times [0,1]$ such that $\partial N=S\cup(F_0\times 0)\cup(-F_1\times 1)$.

Indeed, let $U_k$ be a closed tubular neighborhood of $L_k=\partial F_k$ in $F_k$ 
for $k=0,1$. Let $F'_k$ be the closure of $F_k \setminus U_k$. Deforming if 
necessary $S$, we can assume that $S$ meets $\partial (M\times [0,1])$ precisely 
along $\partial S=(L_1\times 1)\cup(-L_0\times 0)$.  Let 
$U=S\times D^2$ be a closed tubular neighborhood of $S$ in $M\times [0,1]$. Deforming if necessary $U$, we can assume that
$U\cap (F_k \times k)=U_k \times k$ for $k=0,1$. Let $Y$ be the closure of $(M\times [0,1]) \setminus U$. Then 
$Y$ is a compact oriented 4-manifold with boundary and $F'_k\times k\subset\partial Y$ for $k=0,1$.

We define a continuous map $f\colon\partial Y\to S^1$ as follows.
For $k=0,1$, let $F'_k\times[-1,1]$ be a closed tubular neighborhood of 
$F'_k\times k$ in $Y\cap (M\times k)\subset \partial Y$. Then, $f$ restricted
to $F'_k\times[-1,1]$ is given by $f(x,t)=e^{i\pi t}$ for $x\in F'_k, t\in 
[-1,1]$. On $S\times \partial D^2 \subset \partial Y$, the map $f$ is such that 
$f^{-1}(1)=S\times \star$ for some $\star \in \partial D^2$. Finally,
$f(x)=-1$ for all $x\in\partial M\times [0,1]$ and all $x\in((M\times k)\cap Y)\setminus(F'_k\times[-1,1])$ where $k=0,1$. By
elementary obstruction theory, the   map $f\colon\partial Y\to S^1$ 
extends to   $Y$ if and only if there is a
homomorphism $\phi\colon H_1(Y)\to\Zz$ such that $\phi\circ i_*=f_*$, 
where $i$ is the inclusion $\partial Y\hookrightarrow Y$.
Using the exact homology sequence of the pair $(M\times[0,1],Y)$, the 
excision theorem, and the assumption $H_2(M)=0$, we obtain  that $H_3(Y)=0$ and 
$H_2(Y )=\Zz^m$ where $m$ is the 
number of 
components of $L_0$ (and of $L_1$). A basis   of $H_2(Y)$ is given 
by the homology classes of $m$ tori $T_1,\dots, T_m \subset \partial Y$ 
forming $\partial(U\cap (M\times 0))$.  
We have  $H_1(Y,\partial Y)= H^3(Y)=0$ and $H_2(Y,\partial Y)=H^2(Y )=\Zz^m\oplus G$ where $G$ is a finite abelian group. 
The summand $\Zz^m\subset H_2(Y,\partial Y)$ has a basis $y_1,\dots, y_m$ dual to the basis $[T_1],\dots,[T_m]$ of $H_2(Y)$.
The homological  sequence of the pair $(Y,\partial Y)$ yields
$$
H_2(Y,\partial Y)\stackrel{\partial}{\longrightarrow}H_1(\partial Y)\stackrel{i_*}{\longrightarrow}H_1(Y)\longrightarrow 0.
$$
Clearly, $f_*(\partial (G))=0$. Using the assumption $\partial M\neq\emptyset$,
it is easy to construct for each  $j=1,\dots,m$, a loop in 
$f^{-1}(-1)\subset\partial Y$ piercing $T_j$ once and disjoint from the other 
$m-1$ tori. This loop represents $\partial(y_j)$ mod $\partial (G)$.
Therefore, $f_*(\partial (y_j))=0$ for all $j$. Thus, the obstruction to the 
extension of $f$ to $Y$ mentioned above is $0$. Let $\widetilde f\colon Y\to S^1$
be a  continuous extension of $f$.  Deform $\widetilde f$ so that $1$ is one of its regular values.
Then the 3-manifold $N=\widetilde f^{-1}(1)$ satisfies the conditions of Claim 1.

Set $H'=H_1(F'_0;R)\oplus H_1(F'_1;R)$, which we identify  with $H=H_1(F_0;R)\oplus H_1(F_1;R)$ via the inclusion homomorphism.
Let $\widetilde K$ (resp. $K$) be the kernel of the inclusion homomorphism $H_1(\partial N;Q)\to H_1(N;Q)$ (resp. $H\otimes Q\to H_1(N;Q)$),
where $\otimes=\otimes_R$ and $Q=Q(R)$ denotes the field of fractions of $R$. By the standard argument using the  Poincar\'e-Lefschetz duality,
the dimension of $\widetilde K$ is half of the dimension of $H_1(\partial N;Q)$. Furthermore, one easily checks that both the kernel and the cokernel of
the inclusion homomorphism $H\otimes Q\to H_1(\partial N;Q)$ have dimension $m-1$. Therefore,
$$
\dim_Q K\ge\dim_Q\widetilde K=\frac{1}{2}\dim_Q H_1(\partial N;Q)=\frac{1}{2}\dim_Q(H\otimes Q).
$$
We now use this fact to show a second claim. The proof is adapted from \cite[p. 89]{Lick}. 

{\em Claim 2:\/} There is an $R$-basis $x_1,\dots,x_{2g}$ of $H$ such that $x_i$ maps to zero in $H_1(N;Q)$ for all $i>g$.

Observe first that $H$ is a free $R$-module of rang $2g$ where  $g$ is the genus of $\partial N$. Then  $H\otimes Q$ is a vector space over $Q$ of dimension $2g$ and  $\dim_Q K\ge g$. Pick a $g$-dimensional subspace $E$ of $K$. Clearly, $E$ admits a $Q$-basis consisting of elements in $H$:
just take any $Q$-basis of $E$ and multiply its vectors by non-zero scalars. Let $ E_0$ be the $R$-span of these elements in $H$.
Since $R$ is a principal ideal domain, $H/E_0=F\oplus T$  where $F$ is a free $R$-module of rank $g$
and $T$ a torsion $R$-module. Let $\widetilde T $ be the pre-image of $T$ under the projection  $H  \to H/E_0$.  Then
$E_0\subset\widetilde T\subset H\cap E$ and $\widetilde T/E_0= T$. Since $R$ is a principal ideal domain and
$H$ is free, $\widetilde T$ is free as well. Since the sequence $0\to \widetilde T\to H\to F\to 0$ is exact and $F$ is free, a basis for
$\widetilde T$ can be completed  to an $R$-basis of $H $ which satisfies
the conditions of Claim 2.

The lemma now follows from one last claim.

{\em Claim 3:\/} If $a,b\in H$ map to zero in $H_1(N;Q)$, then $\widetilde\psi(a,b)=\vartheta(a,b)=0$.

Indeed, if $a,b\in H$ map to zero in $H_1(N;Q)=H_1(N;R)\otimes Q$, then $r\cdot a$ and $r'\cdot b$ map to zero in $H_1(N;R)$ for some
non-zero $r,r'\in R$. By $R$-bilinearity of $\widetilde\psi$ and $\vartheta$  and the assumption that $R'$ is an integral domain,
it is enough to consider the case where $a,b\in H$ map to zero in $H_1(\partial N;R)$. We have $a=a_0\oplus a_1$ and $b=b_0\oplus b_1$ with 
$a_0, b_0\in H_1(F_0;R)$ and $a_1, b_1 \in  H_1(F_1;R)$.
Consider the following inclusion homomorphisms
$$
H \to H_1(\partial N;R)\to H_1(N;R) \to H_1(M\times [0,1];R)=H_1(M;R)=V .
$$
Clearly, the composition is given by $x_0\oplus x_1\mapsto 
d_0(x_0)+d_1(x_1)$. Since $a,b$ are in the kernel of this composition,
$d_0(a_0)+d_1(a_1)=d_0(b_0)+d_1(b_1)=0$. Hence,
$$
\widetilde\psi(a,b)=-\psi(d_0(a_0),d_0(b_0))+\psi(d_1(a_1),d_1(b_1))=0.
$$
By the assumptions on $a,b$, there are 2-chains $\alpha,\beta$ in $N$ such that 
$\partial\alpha=a_0+a_1$ and $\partial\beta=b_0+b_1$.
Let $B_k$ be a 2-cycle in $M\times k$ such that $\partial B_k=b_k-d_k(b_k)$  for $k=0,1$. Then
$$
\vartheta(a,b)=\vartheta_1(a_1,b_1)-\vartheta_0(a_0,b_0)= a_1^+\cdot_{M\times 1} B_1 - a_0^+ \cdot_{M\times 0} B_0.
$$
The equality $d_0(b_0)+d_1(b_1)=0$ implies that there is 
a 2-chain  $Z$ in $ \partial M\times[0,1]$ such that
$\partial Z=d_0(b_0)+d_1(b_1)$. Since $Z$ is disjoint from $a_0^+$ and 
$a_1^+$, and $a_k^+$ is disjoint from $B_{\ell}$ for $k\neq \ell$,
$$ 
\vartheta(a,b)=(a_0^++a_1^+)\cdot_{\partial(M\times[0,1])} (B_0+B_1+Z).
$$
Here we   used the fact that the orientation on 
$\partial(M\times[0,1])$ matches the one on $M\times 1$ 
and is opposite to the one on $M\times 0$. There is a map $N\to (M\times[0,1])\setminus N$ 
extending the push in the positive normal direction   $F'_k\to (M\times k)\setminus F'_k$ for $k=0,1$.
Let $\alpha^+$ be the image of $\alpha$ under this map. Then 
$$
\vartheta(a,b)=\alpha^+\cdot_{M\times[0,1]}(B_0+B_1+Z)= \alpha^+\cdot_{M\times[0,1]}\beta,
$$
since $B_0+B_1+Z-\beta$ is a $2$-cycle, and therefore a $2$-boundary, in 
$M\times[0,1]$. Finally, $\beta\subset N$ and $\alpha^+\subset (M\times[0,1])\setminus N$ are 
disjoint, so $\alpha^+\cdot_{ M\times[0,1]}\beta=0$.
This concludes the proof.
\end{proof}

The following theorem generalizes the results of Fox-Milnor \cite{F-M} for knots in $S^3$.
\begin{thm}\label{irg}
Let $L_0,L_1$ be concordant homologically trivial links in a quasi-cylinder $(M,V)$ over a principal ideal domain $R$ 
such that $M$ is compact and $H_2(M)=0$. Let $\psi_1,\dots, \psi_n\colon V\times V \to R$ be bilinear forms such that
$\psi_u$ is symmetric for $u=1,\dots, p$ and skew-symmetric for $u=p+1,\dots, n$. Then for some $f\in R[t^{1/2},s_1,\dots, s_n]$,
$$
\Delta_{L_0,\psi_1,\dots,\psi_n}(t^{-1},-s_1t^{-1/2},\dots,-s_p t^{-1/2}, s_{p+1} t^{-1/2},\dots,s_nt^{-1/2})\times
$$
$$
\times\Delta_{L_1,\psi_1,\dots, \psi_n}(t ,s_1t^{ 1/2},\dots, s_nt^{ 1/2})
$$
$$
=f(t^{-1},-s_1,\dots, -s_p, s_{p+1},\dots, s_n)\, f(t,s_1,\dots, s_n).
$$  
\end{thm}
\begin{proof} By Lemma \ref{lemma:conc}, the matrices of the bilinear pairings  
$\vartheta= (-\vartheta_0)\oplus\vartheta_1$ and
$$
\widetilde\psi= -(\sum_u s_u\psi_u \circ(d_0\times d_0))\oplus(\sum_u s_u\psi_u\circ(d_1\times d_1))
$$
with respect to a certain basis of $H=H_1(F_0;R)\oplus H_1(F_1;R)$ have the form
$$
\Theta= \begin{pmatrix}\star&A\cr 
B&0\end{pmatrix}\quad\hbox{and}\quad \widetilde\Psi=\sum_u s_u \widetilde\Psi_u= \begin{pmatrix}\star& \sum_u s_u C_u\cr 
\sum_u s_u C'_u&0\end{pmatrix} ,
$$
where $A,B, C_u, C'_u$  are square matrices over $R$ of  equal size. Note that $C'_u =C^T_u$ for $u=1,\dots, p$ and $C'_u=- C^T_u$ for $u=p+1,\dots, n$.

Let $m$ be the number of components of $L_0$ (and of $L_1$). By Proposition \ref{prop:properties},
$$
\Delta_{L_0, \psi_1,\dots,  \psi_n}(t^{-1},-s_1t^{-1/2},\dots,-s_p t^{-1/2}, s_{p+1} t^{-1/2}, \dots, s_nt^{-1/2})\times 
$$
$$
\times \Delta_{L_1,\psi_1,\dots, \psi_n}(t ,s_1t^{ 1/2},\dots, s_nt^{ 1/2})
$$
\begin{eqnarray*}
&=&(-1)^{m-1}\Delta_{L_0,\psi_1,\dots, \psi_n}(t ,s_1t^{ 1/2},\dots, s_nt^{ 1/2})\,\Delta_{ L_1,\psi_1,\dots, \psi_n}(t ,s_1t^{ 1/2},\dots, s_nt^{ 1/2})\\
&=&\begin{vmatrix}t^{1/2}\Theta-t^{-1/2}\Theta^T+\widetilde\Psi\end{vmatrix}\\
&=&\begin{vmatrix}\star&t^{1/2}A-t^{-1/2}B^T+ \sum_u s_u C_u\cr t^{1/2}B-t^{-1/2}A^T-\sum_u s_u C'_u&0\end{vmatrix}\\
&=&f(t^{-1},-s_1,\dots, -s_p, s_{p+1},\dots, s_n)\, f(t,s_1,\dots, s_n),
\end{eqnarray*}
where $f(t,s_1,\dots, s_n)=\begin{vmatrix}t^{1/2}A-t^{-1/2}B^T+\sum_u s_u C_u\end{vmatrix}$.
(The sign $(-1)^{m-1}$ disappears because of the minuses in the definition of the forms $\vartheta $ and $ \widetilde\psi$.) \end{proof}

\section{Signatures and derived invariants}\label{section:signature}

\subsection{Signatures}
The classical Murasugi-Tristram-Levine signature of a link $L$ in the 3-ball is the function $\sigma_L\colon S^1\to\Zz$
whose value on $\omega\in S^1\subset\Cc$ is the signature of the Hermitian matrix
$(1-\omega)\Theta+(1-\overline\omega)\Theta^T$, where $\Theta$ is a Seifert matrix of $L$. This function is a well-defined invariant of $L$. It is
a concordance invariant away from the roots of $\Delta_L$ on $S^1$. We now extend these results to our setting.
 
Consider a quasi-cylinder $(M,V)$ over $R=\Rr$. Fix $p$ symmetric bilinear forms $\psi_1,\dots, \psi_p\colon V\times V \to \Rr$ and
$ n-p$ skew-symmetric bilinear forms $\psi_{p+1},\dots, \psi_{n}\colon V\times V \to \Rr$. Let $L$ be a homologically trivial link in $M$ and 
$(H, \vartheta, d)$ be the Seifert triple associated with a Seifert surface for $L$.  The {\em signature of $L$\/} is the function
$$
\sigma_{L,\psi_1,\dots, \psi_n}\colon S^1\times \Rr^n\to\Zz
$$
sending a tuple $(\omega\in S^1 , \lambda=(\lambda_1,\dots, \lambda_n) \in\Rr^n)$ to  the signature of the Hermitian form
$$
(1-\omega)\vartheta+(1-\overline\omega)\vartheta^T+(\sum_{u=1}^p \lambda_u\psi_u+ i\sum_{u=p+1}^n\lambda_{u}\psi_u)\circ(d\times d)
$$
on $\Cc \otimes_{\Rr} H$.  Using Theorem \ref{thm1}, one easily checks that $\sigma_L$ does not depend on the choice of the Seifert surface
(see e.g. \cite[Chapter 8]{Lick} for a proof which extends to our setting). Thus, it is a well-defined isotopy invariant of $L$.

\begin{thm}\label{thm:sign}
Let $L_0,L_1$ be concordant homologically trivial links in a quasi-cylinder $(M,V)$ over $\Rr$
such that  $M$ is compact and $H_2(M)=0$. Then
$$
\sigma_{L_0, \psi_1,\dots, \psi_n}(\omega,\lambda)=\sigma_{L_1, \psi_1,\dots, \psi_n}(\omega,\lambda)
$$
for all $\omega \neq 1$ and $\lambda\in\Rr^n$ such that both 
$\Delta_{L_0, \psi_1,\dots, \psi_n}$ and $\Delta_{L_1, \psi_1,\dots, \psi_n}$ do not vanish on 
$\left(\omega, \xi \lambda_1,\dots, \xi \lambda_p, i\xi \lambda_{p+1},\dots,i\xi \lambda_{n}\right)$ where $\xi=(1-\omega^{-1})^{-1}$.
\end{thm}
\begin{proof}  We shall use the notation introduced in the proof of Theorem \ref{irg}.
Clearly, $\sigma_{L_1,\psi_1,\dots, \psi_n}(\omega,\lambda)-\sigma_{L_0,\psi_1,\dots, \psi_n}(\omega,\lambda)=\sgn(\Phi)$, where
 \begin{eqnarray*}
\Phi &=&(1-\omega)\Theta+(1-\overline\omega)\Theta^T+  \sum_{u=1}^p  \lambda_u \Psi_u +  i \sum_{u=p+1}^n   \lambda_{u}   \Psi_u \\
&=&\begin{pmatrix}
\star&(1-\omega)A+(1-\overline\omega)B^T+ C \cr 
(1-\omega)B+(1-\overline\omega)A^T+C'  &0
\end{pmatrix},
\end{eqnarray*}
with $A,B,C, C'$ square matrices over $\Cc$ of equal size. Therefore, $\sgn(\Phi)=0$ unless $\Phi$ is degenerate. We have 
$$
\det 
\Phi=\pm \prod_{k=0,1} \begin{vmatrix}(1-\omega)\Theta_k+(1-\overline\omega)
\Theta_k^T+\sum_{u=1}^p \lambda_u \Psi_{k,u}  +    i\sum_{u=p+1}^n    {\lambda_u} \Psi_{k,u}\end{vmatrix},
$$
where $\Theta_k$ and $\Psi_{k,u}$ are the matrices of the forms $\vartheta_k$ and $\psi_u \circ  (d_k \times d_k)$ on $ H_1(F_k;\Rr)$.
For $k=0,1$, the $k$-th  determinant on the right-hand side is equal to
$$
 \omega^{-r_k/2}
(1-\omega )^{r_k}\Delta_{L_k,   \psi_1,\dots, \psi_n}\left(\omega, \xi \lambda_1,\dots, \xi \lambda_p, i\xi \lambda_{p+1},\dots,i\xi \lambda_{n}
\right),
$$
where $r_k=\dim H_1(F_k;\Rr)$. This proves the theorem.
\end{proof}
  
\subsection{Further invariants}
We assume in this subsection that the ground ring $R$ is a field and $W$ is a vector space over $R$.
More invariants of Seifert triples can be obtained using the following construction.
A Seifert triple $(H,\vartheta,d)$ over $W$ gives a Seifert triple 
$(H',\vartheta',d')$ over any submodule $W'$ of $W$ by $H'=d^{-1} (W')$, 
$\vartheta'=\vartheta\vert_{H'\times H'}$, and $d'=d \vert_{H'}$. The latter triple is 
said to be a {\em restriction\/} of $(H,\vartheta,d)$. Note that equivalent Seifert 
triples   may give non-equivalent restrictions. To handle this, we introduce a 
notion of stable equivalence for Seifert triples.

We say that a  Seifert triple $(H_2,\vartheta_2,d_2)$ over $W$ is 
obtained from a  Seifert triple $(H_1,\vartheta_1,d_1)$ over  $W$ by a {\em trivial enlargement\/} (and
$(H_1,\vartheta_1,d_1)$ is obtained from $(H_2,\vartheta_2,d_2)$ by a  {\em trivial reduction\/}) if 
$H_2= H_1 \oplus Rb$, $d_2|_{H_1} =d_1$, $d_2(b)=0$, $\vartheta_2|_{H_1\times H_1}=\vartheta_1$,
$\vartheta_2(H_1,b)= \vartheta_2(b, H_1)=\vartheta_2(b, b)=0$. Thus, a matrix of $\vartheta_2$ is obtained from a matrix 
of $\vartheta_1$ by adding a zero row and a zero column. Two Seifert triples over $W$ are 
{\em stably equivalent\/} if they can be related by (a finite sequence of) isomorphisms, 
elementary enlargements and reductions, and trivial enlargements and reductions.

It is easy to check that stably equivalent Seifert triples over $W$
restrict  to stably equivalent Seifert triples over submodules of 
$W$. Therefore a  stable equivalence invariant of Seifert triples
generates a family of such invariants by applying it to all possible 
restrictions of a given Seifert triple.

Given a Seifert triple $(H,\vartheta,d)$ over $W$, the associated  
polynomial $\det(t^{1/2}\Theta-t^{-1/2}\Theta^T+ t^{-1/2}\Psi)$ as in Section \ref{sec2.2}
is not preserved under trivial enlargements. The module presented by the matrix 
$t\Theta-\Theta^T+ \Psi$ is preserved up to taking direct sums with free $R'[t,t^{-1}]$-modules 
of finite rank. The sequence of elementary ideals of this module is preserved up 
to shifts of the index. 

The signatures of Seifert triples are easily seen to be invariant under stable equivalence. This generates 
a family of stable equivalence invariants obtained by taking the signatures of the restrictions. 

Applying the constructions above to homologically trivial links in a quasi-cylinder $(M,V)$ over $\Rr$, we obtain 
{\em derived signatures\/} numerated by the subspaces of $V$. They are isotopy invariants. We do not know whether
they are concordance invariants or not.

\section{The multivariable case}\label{section:colored}

The classical theory of Seifert surfaces for oriented links in $S^3$ has been extended to $\mu$-colored links in $S^3$ using `C-complexes'
(see \cite{Cim,C-F,Coop,Coo}). The aim of the present section is to sketch a further extension of this theory to $\mu$-colored links in quasi-cylinders.

\subsection{Colored links}
Let $\mu$ be a fixed positive integer. A {\em $\mu$-colored link\/} $L=L_1\cup\dots\cup L_\mu$ in an oriented $3$-manifold $M$ is
an oriented link in the interior of $M$ together with a surjective map assigning to each component of $L$ a color in $\{1,\dots,\mu\}$.
The sublink $L_i$ is constitued by the components of $L$ with color $i$ for $i=1,\dots,\mu$.
We shall say that two colored links $L,L'$ in $M$ are {\em isotopic\/} if there is an ambient isotopy between $L$ and $L'$, fixing $\partial M$, and
preserving the orientation and color of every component.
A $\mu$-colored link $L=L_1\cup\dots\cup L_\mu$ is {\em homologically trivial\/} if $[L_i]=0$ in $H_1(M)$ for all $i=1,\dots,\mu$.

Note that a $1$-colored link is an ordinary link, as defined in Section \ref{section:lk}. Setting $\mu=1$ in the present section, we obtain  the
theory developed in the previous sections.

\begin{figure}[Htb]
\labellist\small\hair 2.5pt
\pinlabel {$F_i$} at 15 75
\pinlabel {$F_j$} at 380 75
\endlabellist
\centerline{\psfig{file=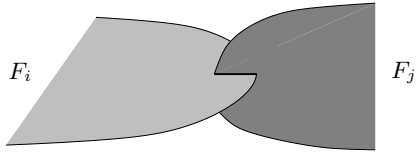,height=2cm}}
\caption{A clasp intersection.}
\label{fig:clasp}
\end{figure}

\subsection{C-complexes}
A {\em C-complex\/} for a $\mu$-colored link $L=L_1\cup\dots\cup L_\mu$ in an oriented $3$-manifold $M$ is a union $F=F_1\cup\dots\cup F_\mu$
of surfaces in $M$ such that $F$ is connected, and the following conditions hold:
\begin{romanlist}
\item{for all $i$, $F_i$ is a Seifert surface for $L_i$;}
\item{for all $i\neq j$, $F_i\cap F_j$ is either empty or a union of clasps (see Figure \ref{fig:clasp});}
\item{for all $i,j,k$ pairwise distinct, $F_i\cap F_j\cap F_k$ is empty.}
\end{romanlist}
In the case $\mu=1$, a C-complex for $L$ is simply a Seifert surface for $L$.

In order to have a C-complex, a $\mu$-colored link clearly needs to be homologically trivial. One easily checks that it is the only obstruction:
every homologically trivial $\mu$-colored link $L=L_1\cup\dots\cup L_\mu$ in an oriented $3$-manifold has a C-complex. 
Indeed, by Proposition \ref{prop:surface}, every sublink $L_i$ admits a Seifert surface $F_i$. Then, by \cite[Lemma 1]{Cim}, each $F_i$ can be isotoped keeping its boundary fixed to give a C-complex for $L$.

\begin{figure}[Hb]
\labellist\small\hair 2.5pt
\pinlabel {$T2$} at 330 200
\pinlabel {$T3$} at 1140 200
\endlabellist
\centerline{\psfig{file=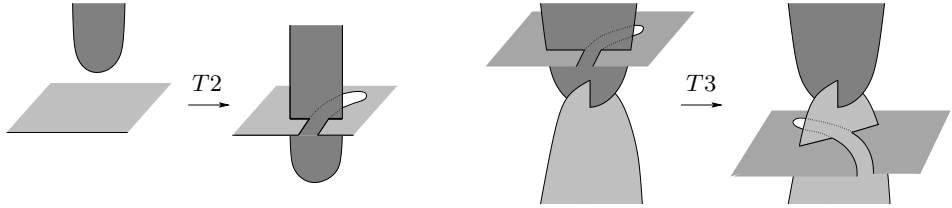,height=2.7cm}}
\caption{The transformations $T2$ and $T3$ in Proposition \ref{prop:Cmoves}.}
\label{fig:Sequiv}
\end{figure}

\begin{prop}\label{prop:Cmoves}
Let $F$ and $F'$ be C-complexes for isotopic colored links in a quasi-cylinder $(M,V)$ over $R$.
If $H_2(M)=0$, then $F$ and $F'$ can be transformed into each other by a finite number of the following operations and their inverses:
\begin{Tlist}
\setcounter{Tlistc}{-1}
\item{Ambient isotopy keeping $\partial M$ fixed;}
\item{surgery on one surface;}
\item{addition of a ribbon intersection, followed by a `push along an arc' through this intersection (see Figure \ref{fig:Sequiv});}
\item{the transformation described in Figure \ref{fig:Sequiv}.}
\end{Tlist}
\end{prop}
\begin{proof}
By the first move, it may be assumed that $\partial F_i=\partial F'_i=L_i$ for all $i$. Since $H_2(M)=0$, $F_i$ and $F_i'$ are related by ambient isotopies (keeping $L_i$ fixed) and surgeries. Clearly, a surgery on $F_i$ can be performed avoiding $F\setminus F_i$, giving move $T1$.
Now, for every ambient isotopy between $F_i$ and $F_i'$, we can apply \cite[Lemma 3]{Cim}, whose proof extends to our setting: such an ambient isotopy can
be induced by a finite sequence of   moves $T0$, $T2$, $T3$ and their inverses.
\end{proof}

\subsection{Seifert forms for colored links}
\begin{figure}[Htb]
\labellist\small\hair 2.5pt
\pinlabel {$F_i$} at 85 160
\pinlabel {$F_j$} at 300 170
\pinlabel {$\eps_i$} at 33 111
\pinlabel {$\eps_j$} at 265 -10
\pinlabel {$F^\eps$} at 465 170
\endlabellist
\centerline{\psfig{file=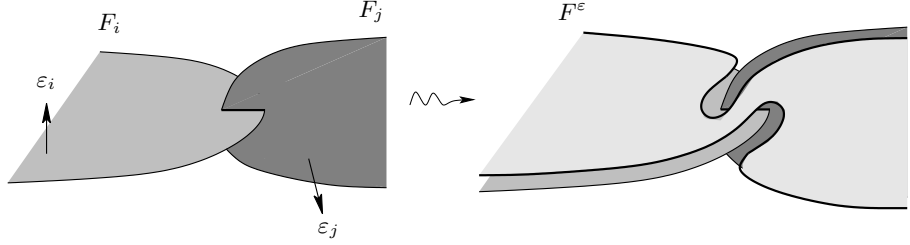,height=3cm}}
\caption{The surface $F^\eps$ near a clasp; the arrow off $F_i$ indicates the $\eps_i$-normal direction on $F_i$ in $M$.}
\label{fig:eps}
\end{figure}

Let us now define the corresponding generalization of the Seifert form. Let as above $R$ be an arbitrary commutative ring with unit.
Let $N_i=F_i\times[-1,1]$ be a bicollar neighborhood of $F_i$ in the interior of $M$. Given a sign $\eps_i=\pm 1$, let $F_i^{\eps_i}$ be the translated
surface $F_i\times\{\eps_i\}\subset N_i$. Also, let $T(L_i)$ be a tubular neighborhood of $L_i$ in $\Int(M)$, and let $Y$ be the complement of
$\bigcup_{i=1}^\mu\Int(N_i\cup T(L_i))$ in $M$. Given a sequence $\eps=(\eps_1,\dots,\eps_\mu)$ of $\pm 1$, set
$$
F^\eps=\bigcup_{i=1}^\mu F_i^{\eps_i}\cap Y.
$$
See Figure \ref{fig:eps} for an illustration of $F^\eps$ near a clasp. Since all the intersections are clasps,
there is an obvious homotopy equivalence between $F$ and $F^\eps$
inducing an isomorphism $H_1(F;R)\to H_1(F^\eps;R)$, $a\mapsto a^\eps$.
Note also that $F^\eps$ is a smooth surface, endowed with a canonical orientation:
the orientation that matches the one on $F_i$ if
and only if $\eps_i=+1$. Hence, we have a well-defined Seifert form $\vartheta_{F^\eps}$ on $H_1(F^\eps;R)$ as in 
Section \ref{section:forms}. Therefore, each choice of signs $\eps=(\eps_1,\dots,\eps_\mu)$ leads to a Seifert form $\vartheta^\eps$ and to an intersection
form $\varphi^\eps$ on $H_1(F;R)$ defined by
$$
\vartheta^\eps(a,b)=\vartheta_{F^\eps}(a^\eps,b^\eps)\quad\hbox{and}\quad\varphi^\eps(a,b)=a^\eps\cdot_{F^\eps}b^\eps
$$
for all $a,b$ in $H_1(F;R)$. These forms are related as follows.

\begin{lemma}\label{lemma:Crelation}
For all $a,b$ in $H_1(F;R)$ and all signs $\eps=(\eps_1,\dots,\eps_\mu)$,
$$
\vartheta^\eps(a,b)-\vartheta^{-\eps}(a,b)=\varphi^\eps(a,b)\quad\hbox{and}\quad\vartheta^\eps(a,b)-\vartheta^{-\eps}(b,a)=d(a)\cdot_{\partial M}d(b),
$$
where $\cdot_{\partial M}$ is the intersection pairing on $\partial M$ and $d\colon H_1(F;R)\to V$ the composition 
of the inclusion homomorphism $H_1(F;R)\to H_1(M;R)$ with the isomorphism $d_V\colon H_1(M;R)\to V$.
\end{lemma}
\begin{proof}
Let $i^\eps\colon H_1(F;R)\to H_1(F^\eps;R)$ denote the isomorphism given by $a\mapsto a^\eps$. As an oriented smooth surface,
$F^\eps$ is diffeomorphic to $-F^{-\eps}$, the surface $F^{-\eps}$ with the opposite orientation.
This leads to a canonical isomorphism $h^\eps\colon H_1(F^\eps;R)\to H_1(F^{-\eps};R)$
such that $h^\eps\circ i^\eps=i^{-\eps}$ and $\vartheta^+_{F^{-\eps}}\circ (h^\eps\times h^\eps)=\vartheta^-_{F^\eps}$.
(Recall that the bilinear form $\vartheta^-_{F^\eps}$ is defined as $\vartheta^+_{F^\eps}=\vartheta_{F^\eps}$ but using $a^-$ instead of $a^+$.)
Therefore:
\begin{eqnarray*}
\vartheta^\eps-\vartheta^{-\eps}&=&\vartheta^+_{F^\eps}\circ(i^\eps\times i^\eps)-\vartheta^+_{F^{-\eps}}\circ(i^{-\eps}\otimes i^{-\eps})\\
	&=&\vartheta^+_{F^\eps}\circ(i^\eps\times i^\eps)-\vartheta^+_{F^{-\eps}}\circ(h^\eps\times h^\eps)\circ(i^{\eps}\times i^{\eps})\\
	&=&(\vartheta^+_{F^\eps}-\vartheta^-_{F^\eps})\circ(i^{\eps}\times i^{\eps}).
\end{eqnarray*}
By formula (2.a) applied to $F^\eps$, this is equal to $\cdot_{F^\eps}\circ(i^{\eps}\times i^{\eps})$ giving the result.
The second equality follows from formula (2.b) in a similar way.
\end{proof}

This result leads to the following definition.
A {\em $\mu$-colored Seifert triple\/} over an $R$-module $W$ is a triple $(H,\{\vartheta^\eps\}_\eps,d)$,
where $H$ is a free $R$-module of finite rank, $\{\vartheta^\eps\}_{\eps}$ a family of $2^{\mu-1}$ bilinear forms on $H$ indexed by
the set
$$
E=\{(\eps_1,\eps_2,\dots,\eps_\mu)\,:\,\eps_1=+1,\,\eps_i=\pm 1\hbox{ for $i>1$}\},
$$
and $d$ a homomorphism $H\to W$. (Note that we don't consider the forms $\vartheta^\eps$ with $\eps_1=-1$ since they can be recovered from the other
forms via Lemma \ref{lemma:Crelation}.)

A $\mu$-colored Seifert triple $(\tilde H,\{\tilde\vartheta^\eps\},\tilde d)$ is obtained from another $\mu$-colored Seifert triple
$(H,\{\vartheta^\eps\},d)$ by a {\em type I elementary enlargement\/} if the following conditions hold:
$\tilde H=H\oplus R a\oplus R b$, $\tilde d|_H=d$, $\tilde d(b)=0$, and there is some index $i$ and some sign $\sigma=\pm 1$ such that
for all $\eps\in E$, the matrix $\Theta^\eps$ for $\vartheta^\eps$ with respect to a basis $h$ of $H$ is related to
the matrix $\tilde\Theta^\eps$ with respect to the basis $h\cup\{a,b\}$ of $\tilde H$ by 
$$
\tilde\Theta^\eps=\begin{pmatrix}\Theta^\eps&\star&0\cr\star&\star&\delta_{\sigma,\eps_i}\cr 0&\delta_{-\sigma,\eps_i}&0\end{pmatrix},
$$
where $\delta$ is the Kronecker symbol. Similarly, one speaks of {\em type II elementary enlargement\/} if the following conditions hold:
$\tilde H=H\oplus R a\oplus R b$, $\tilde d|_H=d$, $\tilde d(b)=0$, and there is some indices $i\neq j$ and some signs $\sigma,\sigma'$ such that
$$
\tilde\Theta^\eps=\begin{pmatrix}
\Theta^\eps&\star&0\cr\star&\star&\delta_{\sigma,\eps_i}\delta_{\sigma',\eps_j}\cr 0&\delta_{-\sigma,\eps_i}\delta_{-\sigma',\eps_j}&0
\end{pmatrix}.
$$
We shall say that two $\mu$-colored Seifert triples over $W$ are {\em equivalent\/} if they can be related by a finite number
of type I and II elementary enlargements (and reductions).

\begin{thm}\label{thmC}
Let $(M,V)$ be a quasi-cylinder over $R$ with $H_2(M)=0$.
For any homologically trivial $\mu$-colored link $L$ in $M$, the equivalence class of the $\mu$-colored Seifert triple of a C-complex for $L$
does not depend on the choice of the C-complex and provides an isotopy invariant of the $\mu$-colored link $L$.
\end{thm}
\begin{proof}
By Proposition \ref{prop:Cmoves}, we are left with the proof that if two C-complexes are related
by transformations $T0$ to $T3$, then the corresponding Seifert triples are equivalent. Obviously, transformation $T0$ does not change the Seifert triple.
It is an easy exercice to check that if a C-complex $\tilde F$ is obtained from a C-complex $F$ via surgery on $F_i$, then the corresponding
Seifert triples are related by a type I elementary enlargement with index $i$. (The sign $\sigma$ is determined by the side of $F_i$ along which
the surgery is performed.) Also, one verifies that transformation $T2$ involving surfaces $F_i$ and $F_j$ corresponds to a type II elementary enlargement
with indices $i,j$, and some signs $\sigma,\sigma'$ given by the orientations of $F_i$ and $F_j$. Finally, consider two C-complexes related by a $T3$
transformation. Then, the two corresponding Seifert triples can be understood as two distinct type II elementary enlargements of some fixed Seifert triple.
This concludes the proof.
\end{proof}

\subsection{The Conway function}
Fix a commutative unital ring $R'$ containing $R$ as a subring, and an $R$-bilinear pairing $\psi\colon V\times V\to R'$.
Consider a homologically trivial $\mu$-colored link $L$ in $M$, and let $(H,\{\vartheta_F^\eps\}_\eps,d)$ be the $\mu$-colored Seifert triple  
associated with a C-complex $F$ for $L$. Let $\Theta_F^\eps$ and $\Psi$ be the matrices of the bilinear forms $\vartheta_F^\eps$ and
$\psi\circ(d\times d)$ with respect to a basis of $H$.  

Let $\Lambda_{R',\mu}$ denote the localization of the ring $R'[t_1^{\pm 1},\dots,t_\mu^{\pm 1}]$ with respect to the multiplicative system generated by
$\{t_i-t_i^{-1}\}_{1\le i\le\mu}$. The {\em (extended) Conway function\/} of $L$ is the element of $\Lambda_{R',\mu}$ defined by
$$
\Omega_{L,\psi}(t_1,\dots,t_\mu)=(-1)^\frac{c-\ell}{2}\prod_{i=1}^\mu(t_i-t_i^{-1})^{\chi(F\setminus F_i)-1}\det(-A_F+\Psi),
$$
where $c$ is the number of clasps in $F$, $\ell=\sum_{i<j}lk_V(L_i,L_j)$, and
$$
A_F=\sum_{\eps\in E}\eps_2\cdots\eps_\mu
\left[t_1t_2^{\eps_2}\cdots t_\mu^{\eps_\mu}\,\Theta_F^\eps+(-1)^\mu(t_1t_2^{\eps_2}\cdots t_\mu^{\eps_\mu})^{-1}(\Theta_F^\eps)^T\right].
$$

\begin{prop}
The extended Conway function is an isotopy invariant of the $\mu$-colored link $L$.
\end{prop}
\begin{proof}
By Proposition \ref{prop:Cmoves} and the proof of Theorem \ref{thmC}, we just need to check that $\Omega_{L,\psi}$ remains unchanged if the C-complex $F$
is transformed via moves $T1$ and $T2$. So, let $\widetilde F$ be a C-complex obtained from $F$ by a surgery on $F_k$. Clearly, the number of clasps $c$
remains the same, while
$$
\chi(\widetilde F\setminus\widetilde F_i)=\begin{cases}\chi(F\setminus F_i)&\text{if $i=k$,}\cr\chi(F\setminus F_i)-2&\text{otherwise.}\end{cases}
$$
Furthermore, the corresponding $\mu$-colored Seifert triples are related by a type I elementary enlargement (with index $i=k$).
Using the equality
$$
\sum_{\eps\in E}\eps_2\cdots\eps_\mu
\left[t_1t_2^{\eps_2}\cdots t_\mu^{\eps_\mu}\,\delta_{\sigma,\eps_k}+(-1)^\mu(t_1t_2^{\eps_2}\cdots t_\mu^{\eps_\mu})^{-1}\delta_{-\sigma,\eps_k}\right]=
$$
$$
=\sum_{\eps_1,\dots,\eps_\mu}\eps_1\cdots\eps_\mu\,t_1^{\eps_1}\cdots t_\mu^{\eps_\mu}\,\delta_{\sigma,\eps_k}
\;=\;\sigma t_k^\sigma\prod_{i\neq k}(t_i-t_i^{-1}),
$$
we get
$$
A_{\widetilde F}=
\begin{pmatrix}
A_F&\star&0\cr
\star&\star&\sigma t_k^\sigma\prod_{i\neq k}(t_i-t_i^{-1})\cr
0&-\sigma t_k^{-\sigma}\prod_{i\neq k}(t_i-t_i^{-1})&0
\end{pmatrix},
\quad
\widetilde\Psi=\begin{pmatrix}
\Psi&\star&0\cr
\star&\star&0\cr
0&0&0
\end{pmatrix}.
$$
Therefore, $\det(-A_{\widetilde F}+\widetilde\Psi)=\prod_{i\neq k}(t_i-t_i^{-1})^2\det(-A_F+\Psi)$. The equality follows.
Now, let $\widetilde F$ be a C-complex obtained from $F$ by a move $T2$ involving $F_k$ and $F_\ell$. The number of clasps $\tilde c$ of $\widetilde F$
is given by $  c+2$, and
$$
\chi(\widetilde F\setminus\widetilde F_i)=\begin{cases}\chi(F\setminus F_i)&\text{if $i=k,\ell$,}\cr\chi(F\setminus F_i)-2&\text{otherwise.}\end{cases}
$$
The corresponding $\mu$-colored Seifert triples are related by a type II elementary enlargement with indices $k,\ell$.
By the equality
$$
\sum_{\eps\in E}\eps_2\cdots\eps_\mu\left[t_1t_2^{\eps_2}\cdots t_\mu^{\eps_\mu}\,\delta_{\sigma,\eps_k}\delta_{\sigma',\eps_\ell}
+(-1)^\mu(t_1t_2^{\eps_2}\cdots t_\mu^{\eps_\mu})^{-1}\delta_{-\sigma,\eps_k}\delta_{-\sigma',\eps_\ell}\right]=
$$
$$
=\sum_{\eps_1,\dots,\eps_\mu}\eps_1\cdots\eps_\mu\,t_1^{\eps_1}\cdots t_\mu^{\eps_\mu}\,\delta_{\sigma,\eps_k}\delta_{\sigma',\eps_\ell}
\;=\;\sigma\sigma' t_k^\sigma t_\ell^{\sigma'}\prod_{i\neq k,\ell}(t_i-t_i^{-1}),
$$
we get
$$
A_{\widetilde F}=
\begin{pmatrix}
A_F&\star&0\cr
\star&\star&\sigma\sigma' t_k^\sigma t_\ell^{\sigma'}\prod_{i\neq k,\ell}(t_i-t_i^{-1})\cr
0&\sigma\sigma' t_k^{-\sigma}t_\ell^{-\sigma'}\prod_{i\neq k,\ell}(t_i-t_i^{-1})&0
\end{pmatrix}.
$$
The invariance follows.
\end{proof}

In the case $\mu=1$, $F$ is a Seifert surface for $L$, and the unique Seifert matrix coincides with the matrix $\Theta$ constructed in Section
\ref{section:forms}. Furthermore, we have $c=\ell=0$, $\chi(F\setminus F_1)=\chi(\emptyset)=0$. Hence, the Conway function is given by
$$
\Omega_{L,\psi}(t_1)=\frac{1}{t_1-t_1^{-1}}\det(-t_1\Theta+t_1^{-1}\Theta^T+\Psi)=\frac{(-1)^{m-1}}{t_1-t_1^{-1}}\Delta_{L,-\psi}(t_1^2),
$$
where $m$ is the number of components of $L$.

If $L'$ is a $\mu$-colored link in an oriented $3$-ball $D^3$ and $L$ is the image of $L'$ 
under an orientation preserving embedding $D^3\hookrightarrow M$, then $\Omega_{L,\psi}(t_1,\dots,t_\mu)=\Omega_{L'}(t_1,\dots,t_\mu)$ is the 
usual Conway function of $L'$, as constructed in \cite{Cim}.

Let us conclude this paragraph with a list of properties of $\Omega_{L,\psi}$ generalizing well-known properties of the Conway function of colored
links in $S^3$. We refer to \cite{Cim} for the proofs which easily extend to our setting.

\begin{prop}$(i)$ Let $L_+$, $L_-$ and $L_0$ be homologically trivial $\mu$-colored links which coincide 
everywhere except in a small $3$-ball where they are related as illustrated below.
(Here, $i$ denotes the color of the strands in the $3$-ball.)

\begin{figure}[h]
\labellist\small\hair 2.5pt
\pinlabel {$i$} at -10 150
\pinlabel {$i$} at 170 150
\pinlabel {$i$} at 245 150
\pinlabel {$i$} at 423 150
\pinlabel {$i$} at 483 150
\pinlabel {$i$} at 655 150
\pinlabel {$L_+$} at 87 -20
\pinlabel {$L_-$} at 340 -20
\pinlabel {$L_0$} at 570 -20
\endlabellist
\centerline{\psfig{file=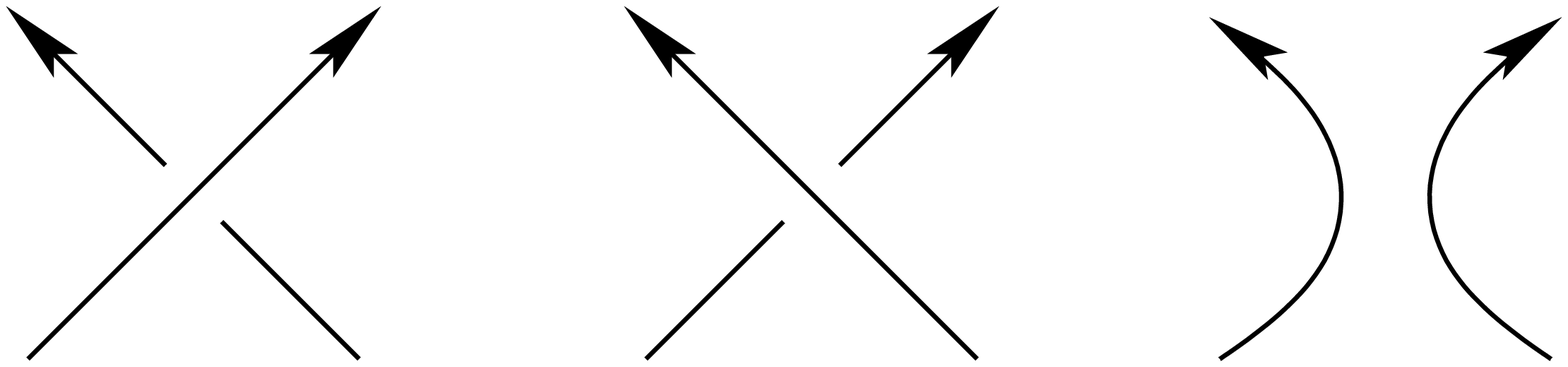,height=1.5cm}}
\end{figure}
\noindent Then, the corresponding Conway functions satisfy the following relation:
$$
\Omega_{L_+,\psi}(t_1,\dots,t_\mu)-\Omega_{L_-,\psi}(t_1,\dots,t_\mu)=(t_i-t_i^{-1})\,\Omega_{L_0,\psi}(t_1,\dots,t_\mu).
$$
$(ii)$ Similarly, if $L_{++}$, $L_{--}$ and $L_{00}$ are homologically trivial $\mu$-colored links which differ by the following local operation,

\begin{figure}[h]
\labellist\small\hair 2.5pt
\pinlabel {$i$} at -20 205
\pinlabel {$j$} at 110 205
\pinlabel {$i$} at 250 205
\pinlabel {$j$} at 383 205
\pinlabel {$i$} at 505 205
\pinlabel {$j$} at 655 205
\pinlabel {$L_{++}$} at 45 -20
\pinlabel {$L_{--}$} at 315 -20
\pinlabel {$L_{00}$} at 575 -20
\endlabellist
\centerline{\psfig{file=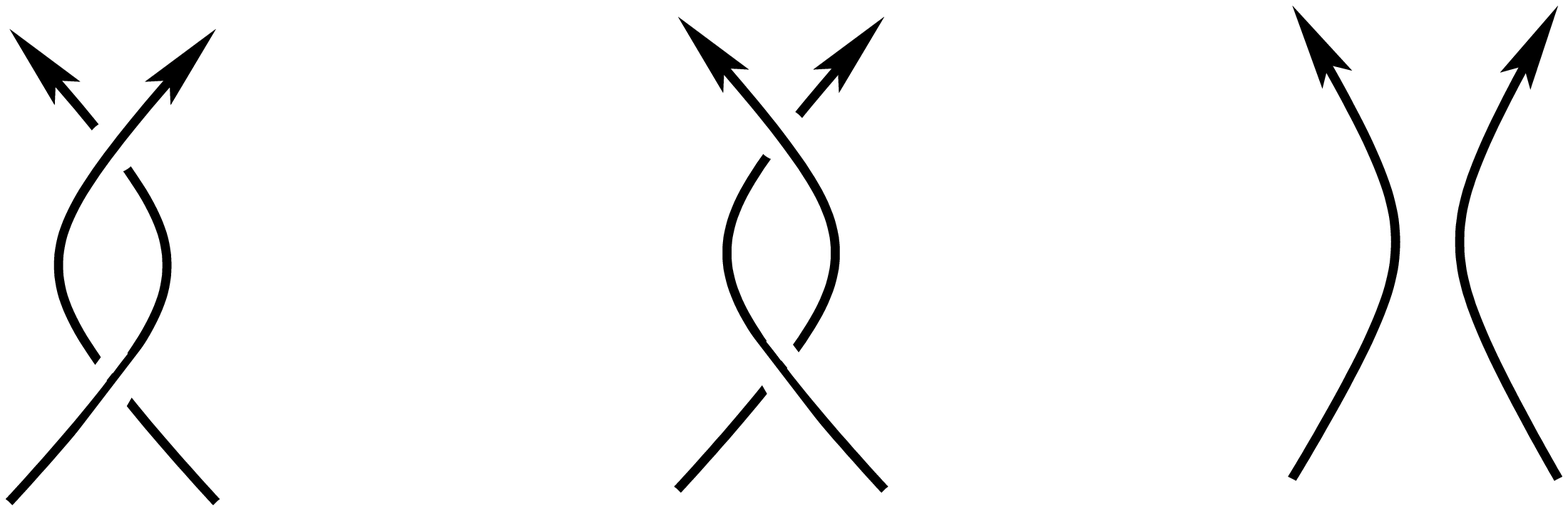,height=1.9cm}}
\end{figure}
then we have the equality
$$
\Omega_{L_{++},\psi}(t_1,\dots,t_\mu)+\Omega_{L_{--},\psi}(t_1,\dots,t_\mu)=(t_it_j+t^{-1}_it^{-1}_j) \, \Omega_{L_{00},\psi}(t_1,\dots,t_\mu).
$$
$(iii)$ For any homologically trivial $\mu$-colored link $L$ with $m$ components,
$$
\Omega_{L,\psi}(t^{-1}_1,\dots,t^{-1}_\mu)=(-1)^m\Omega_{L,\psi'}(t_1,\dots,t_\mu),
$$
where $\psi'$ is the bilinear form given by $\psi'(a,b)=(-1)^\mu\psi(b,a)$.\qed
\end{prop}

\subsection{Multivariable signatures}
As in Section \ref{section:signature}, consider a quasi-cylinder $(M,V)$ over $R=\Rr$, and fix $p$ symmetric bilinear forms
$\psi_1,\dots,\psi_p\colon V\times V\to\Rr$ and $n-p$ skew-symmetric bilinear forms $\psi_{p+1},\dots,\psi_{n}\colon V\times V\to\Rr$.
Let $L$ be a $\mu$-colored homologically trivial link in $M$ and $(H,\{\vartheta^\eps\}_\eps,d)$ be the $\mu$-colored Seifert triple
associated with a C-complex for $L$. Finally, let $T^\mu$ denote the $\mu$-dimensional torus $T^\mu=S^1\times\cdots\times S^1\subset\Cc^\mu$.
The {\em (extended) signature of $L$\/} is the function
$$
\sigma_{L,\psi_1,\dots, \psi_n}\colon T^\mu\times\Rr^n\to\Zz
$$
sending a tuple $(\omega=(\omega_1,\dots,\omega_\mu)\in T^\mu,\lambda=(\lambda_1,\dots, \lambda_n)\in\Rr^n)$ to the signature of the Hermitian form
$$
\sum_{\eps\in E}\Big[(1-\omega_1)\prod_{i>1}(1-\omega^{\eps_i}_i)\vartheta^\eps+
(1-\overline\omega_1)\prod_{i>1}(1-\overline\omega^{\eps_i}_i)(\vartheta^\eps)^T\Big]+\psi
$$
on $\Cc \otimes_{\Rr} H$, where $\psi=(\sum_{u=1}^p \lambda_u\psi_u+ i\sum_{u=p+1}^n\lambda_{u}\psi_u)\circ(d\times d)$.

\begin{prop}
The extended signature is an isotopy invariant of the $\mu$-colored link $L$.
\end{prop}
\begin{proof}
Note that if $\omega_i=1$ for some $i$, then the signature is equal to zero. Therefore, it may be assumed that $\omega_i\neq 1$ for all $i$.
By Theorem \ref{thmC}, we just need to check that the signatures corresponding to equivalent $\mu$-colored Seifert triples are equal.
So, let us assume that a Seifert triple $(\tilde H,\{\tilde\vartheta^\eps\},\tilde d)$ is obtained from another Seifert triple
$(H,\{\vartheta^\eps\},d)$ by a type I elementary enlargement (with index $i=k$).
Using the equality
$$
\sum_{\eps\in E}\Big[(1-\omega_1)\prod_{i>1}(1-\omega^{\eps_i}_i)\delta_{\sigma,\eps_k}+
(1-\overline\omega_1)\prod_{i>1}(1-\overline\omega^{\eps_i}_i)\delta_{-\sigma,\eps_k}\Big]=
$$
$$
=\sum_{\eps_1,\dots,\eps_\mu}\prod_{i=1}^\mu(1-\omega_i^{\eps_i})\,\delta_{\sigma,\eps_k}
=(1-\omega_k^\sigma)\,\prod_{i\neq k}\vert 1-\omega_i\vert^2,
$$
we see that the corresponding Hermitian matrices $\widetilde M$ and $M$ are related by
$$
\widetilde M=
\begin{pmatrix}
M&\star&0\cr
\star&\star&(1-\omega_k^\sigma)\,\prod_{i\neq k}\vert 1-\omega_i\vert^2\cr
0&(1-\overline\omega_k^\sigma)\,\prod_{i\neq k}\vert 1-\omega_i\vert^2&0
\end{pmatrix}.
$$
Since $\omega_i\neq 1$ for all $i$, the signatures of $\widetilde M$ and $M$ coincide by the usual argument.
The invariance of the signature under elementary enlargement of type II follows from the equality
$$
\sum_{\eps\in E}\Big[(1-\omega_1)\prod_{i>1}(1-\omega^{\eps_i}_i)\delta_{\sigma,\eps_k}\delta_{\sigma',\eps_\ell}+
(1-\overline\omega_1)\prod_{i>1}(1-\overline\omega^{\eps_i}_i)\delta_{-\sigma,\eps_k}\delta_{-\sigma',\eps_\ell}\Big]=
$$
$$
=\sum_{\eps_1,\dots,\eps_\mu}\prod_{i=1}^\mu(1-\omega_i^{\eps_i})\,\delta_{\sigma,\eps_k}\delta_{\sigma',\eps_\ell}
=(1-\omega_k^\sigma)(1-\omega_\ell^{\sigma'})\prod_{i\neq k,\ell}\vert 1-\omega_i\vert^2
$$
in the same way.
\end{proof}

In the case $\mu=1$, we obviously get back the extended signatures defined in Section \ref{section:signature}.
If $L'$ is a $\mu$-colored link in an oriented $3$-ball $D^3$ and $L$ is the image of $L'$ 
under an orientation preserving embedding $D^3\hookrightarrow M$, then
$\sigma_{L,\psi}(\omega,\lambda)=\sigma_{L'}(\omega)$ is the multivariable signature of the $\mu$-colored link $L'$, as constructed in \cite{C-F}.

We don't know to which extent the concordance properties of these two special cases (see Theorem \ref{thm:sign} and \cite[Section 7]{C-F})
hold in the general case considered here.

\section{Generalizations}\label{section:gene}

Our invariants of links are defined  under rather strong assumptions: the links 
are supposed to be 
homologically trivial; the ambient manifold, $M$, is supposed to have trivial 
2-homology
and the inclusion homomorphism $H_1(\partial M;R)\to H_1(M;R)$ is supposed to be 
surjective and to 
have a section.  We explain how to weaken these conditions.

\subsection{Homologically non-trivial links}
Let $(M,V)$ be a quasi-cylinder over $R$ with $H_2(M)=0$. Let $h\in H_1(M)$ belong to the image 
of the inclusion homomorphism $H_1(\partial M)\to H_1(M)$. To construct 
non-trivial invariants of links in $M$ representing   $h$, one can proceed as follows. 
Pick a link $L_\ast$ in a cylinder neighborhood $U\subset M$ of $\partial M$ such that 
$[L_\ast]=-h$. Any link $L\subset M$ may be isotopically deformed in $M-U$   uniquely 
up to isotopy in $M-U$. If $L\subset M-U$ and  $[L]=h $, then $\widetilde L=L\cup L_\ast$ is a 
homologically trivial link in $M$. 
The isotopy type of $\widetilde L$ is entirely determined by the isotopy type of $L$  
and the isotopy type of $L_\ast$ in $U$. The invariants of homologically trivial links in $M$ 
defined above may be applied to $\widetilde L$. This yields isotopy invariants of $L$ depending 
on $V$ and $L_\ast$.  In particular, concordance invariants of homologically trivial 
links yield concordance invariants of $L$. Indeed, if two links $L_0,L_1$ in $M$ are 
concordant, then $\widetilde L_0$ and $\widetilde L_1$ are concordant.

\subsection{Generalized quasi-cylinders}
A {\em generalized quasi-cylinder\/} over $R$ is a pair consisting of an oriented 3-manifold $M$ and a submodule $V$ of  
$H_1(\partial M; R)$ such that the inclusion homomorphism $i\colon V\to H_1(M;R)$ is injective. The theory of Seifert triples
associated with surfaces  in quasi-cylinders extend to generalized quasi-cylinders as 
follows. Given an oriented surface $F$ in the interior of $M$, set
$H=j^{-1}(i(V))\subset H_1(F;R)$ where $j$ is the inclusion homomorphism $H_1(F;R)\to H_1(M;R)$. For 1-cycles $a,b$ 
on $F$ representing homology classes $[a], [b] \in H $, set $\vartheta ([a],[b])=lk_V(a^+,b)$. 
This yields a well-defined bilinear form  $\vartheta\colon H\times H\to R$. Applying this construction
to the Seifert surface for a link $L$ in $M$, 
we obtain the Seifert triple $(H,\vartheta, d\colon H\to V)$ of $L$. If $H_2(M)=0$ and $R$ is a field, then the stable 
equivalence class of $(H,\vartheta,d)$ does not depend on the choice of $F$ and 
yields an isotopy invariant of $L$.

\subsection{High-dimensional generalizations}
The constructions of this paper can be easily generalized to  codimension 1 submanifolds of  
odd-dimensional manifolds with boundary and to codimension 2 links in such manifolds.

\subsection{The case of non-connected boundary}
The definitions of linking numbers and generalized Seifert forms given in Sections \ref{section:lk} and
\ref{section:forms} make perfect sense whether $H_2(M)$ is trivial or not (that is, whether $\partial M$
is connected or not). However, the triviality of $H_2(M)$ is needed for Theorem \ref{thm1} to hold. Indeed,
this result is based on the fact that two Seifert surfaces for a link in $M$ can be related by surgeries.
This is clearly not true if $H_2(M)\neq 0$. Therefore, the general theory of Sections
\ref{section:Alex} to \ref{section:colored} does not hold if the boundary of $M$ is non-connected, and
it is very unlikely that any Seifert type invariant can be constructed in this general setting.

Nevertheless, parts of the theory can be developed in the following special case. Let $(M,V)$ be a quasi-cylinder over $R$,
and let us assume that $M$ has exactly two boundary components $\sur$ and $\sur'$, with $V=H_1(\sur;R)$. This is a
natural class of quasi-cylinders, as it contains the prototypical example $M=\sur\times [0,1]$ with $\sur$ closed.
Let $F$ be a Seifert surface in such a quasi-cylinder $(M,V)$, and let $\widetilde\sur$ 
denote a parallel copy of $\sur$ obtained by pushing $\Sigma$ in $\Int(M)\setminus F$. Suppose that there is a solid cylinder $[0,1]\times D^2$ in the interior of $M$ such that
$([0,1]\times D^2)\cap F=\{0\}\times D^2$ and $([0,1]\times D^2)\cap\widetilde\sur=\{1\}\times D^2$. Then  we shall say 
that the surface
$$
F'=(F\setminus (\{0\}\times D^2))\cup([0,1]\times\partial D^2)\cup(\widetilde\sur\setminus (\{1\}\times D^2))
$$
is obtained from $F$ by {\em adding  $\widetilde\sur$ along the arc $[0,1]\times\{0\}$\/}. Here, the orientation of
$\widetilde\sur$ is chosen so that the orientation of $F$ extends to $F'$.

\begin{prop}\label{prop:closed}
Let $(M,V)$ be a compact quasi-cylinder over $R$ with $\partial M=\sur\sqcup\sur'$ and $V=H_1(\sur;R)$.
Any two Seifert surfaces $F, F'$ for isotopic links in a $(M,V)$ can be related by a finite number of ambient isotopies keeping $\partial M$ fixed,
surgeries, and additions of parallel copies of $\sur$ along embedded arcs in $\Int(M)$.
\end{prop}
\begin{proof}
Consider a path $\gamma\colon[0,1]\to M$ such that $\gamma([0,1])\cap\sur=\gamma(0)$, $\gamma([0,1])\cap\sur'=\gamma(1)$, 
and such that $\gamma$ intersects $\sur$, $\sur'$, $F$ and $F'$ transversally. Let us 
assume that $F$ intersects $\gamma$ in $n$ points. Let $\widetilde\sur$ be a
parallel copy of $\sur$ pushed into $M$, disjoint from $F$, and which 
intersects $\gamma$ transversally in $\gamma(t_0)$. Let $t_1$ be the
smallest number such that $\gamma(t_1)\in F$. Consider the surface 
$F_1$ obtained from $F$ by adding  $\widetilde\sur$ along the arc $\gamma([t_0,t_1])$.
Clearly, $F_1$ intersects $\gamma$ in $n-1$ points. Iterating this 
construction, we obtain    a Seifert surface $\widehat F $ for $L$  
disjoint from $\gamma$. Similarly, we obtain a Seifert surface 
$\widehat F'$ from $F'$ disjoint from $\gamma$. Now, consider the 
compact manifold $\widehat M$ given by the complement in $M$ of an open tubular 
neighborhood of $\gamma$. Also, let $\widehat\sur$ be the surface with 
boundary given by $\widehat\sur=\sur\cap\widehat M$. By excision, $H_*(M,\widehat M)=H_*(D^2,S^1)$,
so the homological sequence of $(M,\widehat M)$ reads
$$
0\to H_2(\widehat M)\to H_2(M)\stackrel{\partial}{\to}\Zz.
$$
Since $\partial M$ has exactly two components, one of which is $\sur$,
the inclusion homomorphism $H_2(\sur)\stackrel{i_*}{\to}H_2(M)$ is an 
isomorphism, as well as the composition
$H_2(\sur)\stackrel{\partial\circ i_*}{\longrightarrow}\Zz$. Therefore, 
$\partial$ is an isomorphism, and $H_2(\widehat M)=0$. So, we have two 
Seifert surfaces
$\widehat F$ and $\widehat F'$ in $\widehat M$ for a fixed link $L$ in 
$(\widehat M,\widehat\sur)$, with $H_2(\widehat M)=0$.
By the standard argument,
$\widehat F$ and $\widehat F'$ are related by surgeries in $\Int(\widehat M)\subset\Int(M)$ and by isotopies of $\widehat M$ keeping its boundary fixed.
Such an isotopy obviously extends to an isotopy of $M$ fixing $\partial M$. This concludes the proof.
\end{proof}

Note that $V=H_1(\Sigma;R)$ is endowed with a natural $R$-bilinear form: the intersection form
on $\Sigma$. This leads to the following definition.

Let $W$ be a free $R$-module of finite rank   equipped with   bilinear form $\varphi:W \times W\to R$.
Let $(H,\vartheta,d)$ and $(H',\vartheta',d')$ be two Seifert triples 
over $W$. We shall say that $(H',\vartheta',d')$ is obtained from
$(H,\vartheta,d)$ by a {\em $\varphi$-enlargement\/} (and $(H,\vartheta,d)$ from 
$(H',\vartheta',d')$ by a {\em $\varphi$-reduction\/}) if the following conditions hold:
$H'=H\oplus W$, $d'|_H=d$, $d'|_W=id_W$, $\vartheta'|_{H\times H}=\vartheta$, $\vartheta'|_{H\times W}=0$,
$\vartheta'|_{W\times H}=\varphi\circ(id_W\times d)$ and $\vartheta'|_{W\times W}=0$ or $\varphi$. If $h$ is a basis of $H$
and $w$ a basis of $W$, then $h\cup w$ is a basis of $H'$ and
the matrix $\Theta'$ for $\vartheta'$ with respect to $h\cup w$ is computed from
the matrix $\Theta$ for $\vartheta$ with respect to $h$ by
$$
\Theta'=\begin{pmatrix}\Theta&0\cr 
C&D\end{pmatrix}\quad\hbox{or}\quad\begin{pmatrix}\Theta&0\cr 
C&0\end{pmatrix},
$$
where $C$ is the matrix of $\varphi\circ(id_W\times d)$, and $D$ the matrix of $\varphi$.
We shall say that two Seifert triples over $W$ are {\em $\varphi$-equivalent\/} if they can be related by a finite number of 
isomorphisms, elementary enlargements, elementary reductions, $\varphi$-enlargements and $\varphi$-reductions.

\begin{thm}\label{thm:closed}
Let $(M, V)$ be a quasi-cylinder over $R$ and let us assume that $M$ has exactly two boundary components $\sur$ and $\sur'$,
with $V=H_1(\sur;R)$. Finally, let $\varphi$ denote the intersection form on $V$.
For any homologically trivial link $L\subset M$, the $\varphi$-equivalence class of 
the Seifert triple of a Seifert surface for $L$  does not depend on the 
choice of the surface and provides an isotopy invariant of $L$.
\end{thm}
\begin{proof}
By Proposition \ref{prop:closed}, we just need to check that the addition of a parallel copy of $\sur$ 
induces a $\varphi$-enlargement of the corresponding Seifert triple.
Let $F'$ denote the Seifert surface obtained from $F$ by the addition of 
$\widetilde\sur$ along an arc, and let $\vartheta'$ denote the corresponding 
form. Clearly, $H_1(F')=H_1(F)\oplus H_1(\sur)$, $d'|_{H_1(F)}=d$, $d'|_{H_1(\sur)}=id_{H_1(\sur)}$ and $\vartheta'$ 
restricted to $H_1(F)\times H_1(F)$ is equal to $\vartheta$. Furthermore,
$\vartheta'(a,b)=a^+\cdot_{\partial M}B=0$ for $(a,b)$ in $H_1(F)\times H_1(\sur)$, since $B$ can be chosen to be a thin annulus
$b\times [0,\eta]$ disjoint from $a^+$. For $a,b$ in $H_1(\sur)$,
$$
\vartheta'(a,b)=a^+\cdot_{\partial M}(b\times[0,\eta])= a\cdot_{\partial M}b
$$
if the orientation of $\widetilde\sur$ is induced by the one of $\sur$ and 
$$
\vartheta'(a,b)=a^+\cdot_{\partial M}(b\times[0,\eta])= 0
$$
if the orientation of $\widetilde\sur$ is opposite to the one induced from $\sur$. Finally, for $(a,b)$ in 
$H_1(\sur)\times H_1(F)$, Lemma \ref{lemma:relation} and the above computation give
$$
\vartheta'(a,b)=\overbrace{\vartheta'(b,a)}^{=0}+d'(a)\cdot_{\partial M}d'(b)+\overbrace{a\cdot_{F'} b}^{=0}=a\cdot_{\partial M}d(b).
$$
This concludes the proof.
\end{proof}

Using this theorem, let us now see to which extent the results of Sections \ref{section:Alex} to
\ref{section:colored} hold true in the case under study.

The $R'[t,t^{-1}]$-module $\mathcal{A}_\psi(L)$ is no longer an invariant of $L$ in general.
However, it is an invariant in the special case $R'=R$ and $\psi=-\varphi$, where $\varphi$ is the intersection form on $V$.
Indeed, if $(H',\vartheta',d')$
is obtained from $(H,\vartheta,d)$ by a $\varphi$-enlargement, then the corresponding matrices $\Gamma'=t\Theta'-(\Theta')^T+\Psi'$ and
$\Gamma=t\Theta-\Theta^T+\Psi$ are related by
$$
\Gamma'=\begin{pmatrix}\Gamma&0\cr(t-1)C&tD\end{pmatrix}\quad\hbox{or}\quad \Gamma'=\begin{pmatrix}\Gamma&0\cr(t-1)C&-D\end{pmatrix}.
$$
Since $D$ is congruent to the matrix $\begin{pmatrix}0&1\cr -1&0\end{pmatrix}^{\oplus g}$, $\mathcal{A}_{-\varphi}(L)$
is an invariant of the link $L$. 
Now, consider the element of $R[t^{1/2},t^{-1/2}]$ given by
$$
\widetilde\Delta_L(t)=\Delta_{L,-\varphi}(t)=\det(t^{1/2}\Theta-t^{-1/2}\Theta^T+t^{-1/2}\Psi).
$$
It is   well-defined up to multiplication by $t^g$, where $g$ denotes the genus of $\sur$. Indeed,
if $(H',\vartheta',d')$ is obtained from $(H,\vartheta,d)$ by a $\varphi$-enlargement, then
$$
\det(t^{1/2}\Theta'-t^{-1/2}(\Theta')^T+t^{-1/2}\Psi')=\det(t^{1/2}\Theta-t^{-1/2}\Theta^T+t^{-1/2}\Psi)\cdot\det(\pm t^{-1/2}D).
$$
Since $D$ is a matrix of the intersection form on $\sur$, $\det(\pm t^{-1/2}D)=t^{-g}$, giving the result.
One easily checks the following properties: If $m$ is odd, then $\widetilde\Delta_L(t)\in R[t,t^{-1}]$.
If $m$ is even, then $t^{1/2}\widetilde\Delta_L(t)\in R[t,t^{-1}]$. Finally, $\widetilde\Delta_L(1)=1$ if $L$ is a knot,
and $\widetilde\Delta_L(1)=0$ else.

Proposition \ref{prop:genus1} translates into the inequality
$$
g(L)\ge\frac{1}{2}(\spa\widetilde\Delta_L(t)+1-m).
$$
Furthermore, the Seifert algorithm and Proposition \ref{prop:genus2} extend verbatim to our case.

Generally speaking, the signatures introduced in Section \ref{section:signature} are not invariant under $\varphi$-enlargements.

\subsection*{Acknowledgments}
The first author wishes to thank the UC Berkeley Department of Mathematics for hospitality.
He also expresses his thanks to Mathieu Baillif. The second named author thanks Research Institute for Mathematical Sciences (RIMS, Kyoto) for hospitality during the preparation of this paper.

\end{document}